# SMASH: STRUCTURED MATRIX APPROXIMATION BY SEPARATION AND HIERARCHY [*]

DIFENG CAI[†], EDMOND CHOW[‡], YOUSEF SAAD[§], AND YUANZHE XI[§]

**Abstract.** This paper presents an efficient method to perform Structured Matrix Approximation by Separation and Hierarchy (SMASH), when the original dense matrix is associated with a kernel function. Given points in a domain, a tree structure is first constructed based on an adaptive partitioning of the computational domain to facilitate subsequent approximation procedures. In contrast to existing schemes based on either analytic or purely algebraic approximations, SMASH takes advantage of both approaches and greatly improves the efficiency. The algorithm follows a bottom-up traversal of the tree and is able to perform the operations associated with each node on the same level in parallel. A strong rank-revealing factorization is applied to the initial analytic approximation in the *separation* regime so that a special structure is incorporated into the final nested bases. As a consequence, the storage is significantly reduced on one hand and a hierarchy of the original grid is constructed on the other hand. Due to this hierarchy, nested bases at upper levels can be computed in a similar way as the leaf level operations but on coarser grids. The main advantages of SMASH include its simplicity of implementation, its flexibility to construct various hierarchical rank structures and a low storage cost. Rigorous error analysis and complexity analysis are conducted to show that this scheme is fast and stable. The efficiency and robustness of SMASH are demonstrated through various test problems arising from integral equations, structured matrices, etc.

**1. Introduction.** The invention of the Fast Multipole Method (FMM) [50, 27] opened a new chapter in scientific computing methodology by unraveling a set of effective techniques revolving around the powerful principle of divide-and-conquer. When sets of points are far apart from each other, the physical equations that couple them can be approximately expressed by means of a low rank matrix. Among the many variations to this elegant idea, a few schemes have been developed to gain further efficiency by building hierarchical bases in order to expand the various low-rank couplings. The resulting hierarchical rank structured matrices [12, 14, 22, 31], culminated in $\mathcal{H}^2$ matrices[34, 31], provide efficient solution techniques for structured linear systems (Toeplitz, Hankel, etc.)[54, 56, 60], integral equations [4, 12, 14, 23, 37, 40, 42, 43], partial differential equations [14, 15, 31, 58], matrix equations [14, 25, 26] and eigenvalue problems [8, 55]. Though these methods come under various representations, they all start with a block partitioning of the coefficient matrix and approximate certain blocks with low-rank matrices. The refinements of these techniques embodied in the $\mathcal{H}^2$ [14, 31, 32] and HSS [22, 53, 59] matrix representations take advantage of the relationships between different (numerically) low-rank blocks and use *nested bases* [31] to minimize computational costs and storage requirements. What is often not emphasized in the literature, is that this additional efficiency in the solution phase is achieved at a rather high cost in the construction phase.

Both HSS and $\mathcal{H}^2$ matrices employ just two key ingredients: low-rank approximations and nested bases. The low-rank approximation, or compression, methods exploited in these techniques can be classified into three categories. The first category involves methods that rely on algebraic compression, such as the SVD and the rank–revealing QR (RRQR) [30] which are among the most common approaches. Utilizing these techniques to compress low-rank blocks [10, 32, 59] will result in nested bases that have orthogonal columns and an optimal rank. However, these methods will require the matrix entries to be explicitly available and usually lead to quadratic construction cost [32, 59]. Other compression techniques, such as adaptive cross approximation (ACA) [4, 6, 7], extract a low-rank factorization based only on part of the matrix entries and this leads to a nearly linear construction cost. However, ACA may fail for general kernel functions and complex geometries due to the heuristic nature of the method [13]. To the best of our knowledge, no algebraic approach is able to achieve linear cost for an $\mathcal{H}^2$ or HSS construction with guaranteed accuracy. The methods in the second category rely on information on the kernel to perform the compression. They include methods based on interpolation [14, 33], Taylor expansion [34] or multipole expansion (as in FMM [27, 50]), etc. Although these methods lead to a linear construction cost, they usually yield nested bases whose rank is much larger than the optimal one [35]. Moreover, since bases are stored as dense matrices, these methods suffer from high storage costs [12]. The methods in the third category either combine algebraic compression with the analytic kernel information to take advantage of both, or use other techniques like equivalent densities or randomized methods to obtain a low-rank factorization. For example, hybrid cross approximation (HCA) [13] technique improves the robustness of ACA by applying it only on a small matrix arising from the interpolation of the kernel function. The kernel independent fast multipole methods [1, 61] use equivalent densities to avoid explicit kernel expansions but it is only valid for certain kernels arising in potential theory. The randomized construction algorithms [41, 44, 51, 55] compute the hierarchical rank structured matrices by applying SVD/RRQR to the product of the original matrix and a random matrix and are effective when a fast matrix vector multiplication routine is available.

---

[*]The research of Edmond Chow was supported by NSF under grant ACI–1306573. The research of Yousef Saad and Yuanzhe Xi were both supported by NSF under grants DMS–1216366 and DMS–1521573 and by the Minnesota Supercomputing Institute.

[†]Address: Department of Mathematics, Purdue University, West Lafayette, USA, cai92@purdue.edu
[‡]Address: School of Computational Science and Engineering, Georgia Institute of Technology, Atlanta, USA, echow@cc.gatech.edu
[§]Address: Department of Computer Science & Engineering, University of Minnesota, Twin Cities, USA, {saad,yxi}@umn.edu

**1.1. Contributions.** The aim of this paper is to introduce an efficient and unified framework to construct an $n \times n$ $\mathcal{H}^2$ or HSS matrix based on structured matrix approximation by separation and hierarchy (SMASH). In terms of the three categories discussed above, SMASH belongs to the third category in that it starts with an initial analytic approximation in the *Separation* regime, then algebraic techniques are employed to postprocess the approximation in order to build up a self-similar *Hierarchy*. The main features of SMASH are as follows.

*1. Fast and stable $O(n)$ construction.* SMASH starts with an adaptive partitioning of the computational domain and then constructs a tree structure to facilitate subsequent operations as in [3, 12, 14, 20]. The construction process follows a bottom-up traversal of the tree and is able to compute the bases associated with each node on the same level in parallel. In fact, the construction procedure is entirely local in the sense that the compression for a parent node only depends on the information passed from its children. By combining the analytic compression technique with strong RRQR [30], a special structure is incorporated into the final nested bases. In contrast to the methods used in [44, 61], SMASH is able to set the approximation accuracy to any tolerance. In addition, the nested bases at each non-leaf level can be computed directly in a similar way as the leaf-level operations but on a coarser grid extracted from previous level of compression. Therefore, SMASH is also advantageous relative to one based on the HCA method [11] since it does not need to construct an $\mathcal{H}$ matrix first and then use a recompression scheme to recover the nested bases. SMASH can be easily adapted to construct either an $n \times n$ $\mathcal{H}^2$ or HSS matrix depending on the properties of the underlying applications with $O(n)$ complexity. Equipped with a rigorous error analysis, the guaranteed accuracy/robustness of SMASH is justified by various test examples with complicated geometries (Section 6).

*2. Low storage requirement.* Construction algorithms that use analytic approximations usually lead to high storage costs. SMASH alleviates this issue in several ways. First, instead of storing nested bases as dense matrices, only one vector and one smaller dense matrix need to be saved for each basis. Second, each coupling matrix [12] is a submatrix of the original matrix in this scheme. Therefore, it suffices to store row and column indices associated with the submatrix instead of the whole submatrix explicitly. Finally, the use of strong RRQR [30] can automatically reduce the rank of the nested bases if their columns obtained from the analytic approximation are not numerically linearly independent.

*3. Simplicity and flexibility for approximation of variable order.* Unlike analytic approaches (e.g., FMM) in which farfield approximations and transfer matrices are obtained differently and extra information is needed to compute transfer matrices (cf.[52]), SMASH only requires a farfield approximation, which can be readily obtained for almost all kernels, for example, via interpolation [33]. Moreover, the approximation rank in the compression on upper levels is independent of the rank used in lower levels, which means that approximation rank can be chosen arbitrarily in the compression at any level while still maintaining the $\mathcal{H}^2$ or HSS structure. This is due to the fact that in each level of compression, SMASH produces transfer matrices directly, which is an advantage of algebraic approaches. For interpolation-based constructions, there are restrictions on the maximal degree of basis polynomials in each level in order to maintain the $\mathcal{H}^2$ structure. [32].

**1.2. Outline and Notation.** The paper is organized as follows. In Section 2 we review low-rank approximations ([5, 12, 31]) associated with some kernel functions. Section 3 introduces SMASH for the construction of hierarchical rank structured matrices with nested bases. The approximation error and complexity of SMASH are analyzed in Section 4 and Section 5, respectively. Numerical examples are provided in Section 6 and final concluding remarks are drawn in Section 7.

Throughout the paper, the following notation is used:
- $A$: a dense matrix associated with a kernel function $\kappa$;
- $\hat{A}$: $\mathcal{H}^2$ or HSS approximation to $A$;
- $i = 1 : n$ denotes an enumeration of index $i$ from 1 to $n$;
- $|\cdot|$ denotes the cardinality of a finite set if the input is a set or a vector;
- $\|\cdot\|, \|\cdot\|_F$ denote the $L^2$ norm, Frobenius norm, respectively, and $\|A\|_{\max}$ denotes the elementwise sup-norm of a matrix, i.e.,

$$\|A\|_{\max} := \max_{i,j}|a_{i,j}|, \quad A = [a_{i,j}]_{i,j};$$

- diag(...) denotes a block diagonal matrix;
- Given a tree $\mathcal{T}$, children(i) and lv(i) represent the children and level of node $i$, respectively, where root node is at level 1. The location of a node $i$ at level $l$ is denoted as $l_i$ when enumerated from left to right;
- Let $X$ and $Y$ be two nonempty finite sets of points and $A$ be a matrix whose $(i,j)$th entry is determined by the $i$th point in $X$ and $j$th point in $Y$. If **i** denotes the index set corresponding to a subset $X_\mathbf{i}$ of $X$, then $A|_\mathbf{i}$ denotes the submatrix of $A$ with rows determined by $X_\mathbf{i}$. Furthermore, if index set **j** corresponds to a subset $Y_\mathbf{j}$ of $Y$, then $A|_{\mathbf{i}\times\mathbf{j}}$ denotes a submatrix of $A$ whose rows and columns are determined by $X_\mathbf{i}$ and $Y_\mathbf{j}$, respectively.

**2. Degenerate and low-rank approximations.** Hierarchical rank structured matrices are often used to approximate matrices after a block partitioning such that most blocks display certain (numerical) low-rank characteristics. For matrices derived from kernel functions, a low-rank approximation can be determined when the kernel function can be locally approximated by *degenerate functions* [5]. In this section, we first review this property. For pedagogical reasons, we focus on the kernel function $1/(x-y)$ but more general kernel functions can be handled in a similar way as demonstrated in the numerical experiments section (Section 6).

**2.1. Degenerate expansion.** Consider the kernel function $\kappa(x,y)$ on $\mathbb{C} \times \mathbb{C}$ defined by

$$(2.1) \qquad \kappa(x,y) = \begin{cases} \frac{1}{x-y}, & \text{if } x \neq y, \\ d_x, & \text{if } x = y, \end{cases}$$

where the number $d_x \in \mathbb{C}$ can be arbitrary. If $x$ and $y$ are far from each other (See Definition 2.1 below), then $\kappa(x, y)$ can be well approximated by a degenerate expansion

$$\kappa(x,y) \approx \sum_{k=0}^{r-1} \sum_{l=0}^{k} c_{k,l} \phi_k(x) \psi_l(y),$$

where $\phi_k$ and $\psi_l$ are univariate functions. In fact, interpolation in the $x$ variable yields the simplest, yet most general, way to obtain a degenerate approximation:

(2.2) $$\kappa(x,y) \approx \sum_{k=1}^{r} p_k(x) \kappa(x_k, y),$$

where $x_k$'s are the interpolation points and the $p_k$'s are the associated Lagrange polynomials.

Several ways to quantify the distance between two sets of points that are away from each other have been defined [12, 31, 52]. One of these ([52]), given below, is often referred to. For a bounded nonempty set $S$ of $\mathbb{C}$, let $\delta = \min_{c \in \mathbb{C}} \sup_{s \in S} |s - c|$. Then we refer to the mininizer $c_*$ as the center of $S$ and to the corresponding minimum value $\delta$ as its radius.

DEFINITION 2.1. *Let $X$ and $Y$ be two nonempty bounded sets in $\mathbb{C}$. Let $a \in \mathbb{C}$ and $\delta_a > 0$ be the center and radius of $X$ with $|x - a| \leq \delta_a$, $\forall x \in X$. Analogously, let $b \in \mathbb{C}$ and $\delta_b > 0$ denote the center and radius of $Y$. Given a number $\tau \in (0, 1)$, we say that $X$ and $Y$ are well-separated with separation ratio $\tau$ if the following condition holds*

(2.3) $$\delta_a + \delta_b \leq \tau |a - b|.$$

Fig. 2.1 illustrates two well-separated intervals (centered at $a = 0.5$, $b = 2.5$, respectively) with separation ratio $\tau = 0.5$. Given two sets $X$ and $Y$, if (2.3) only holds for $\tau \approx 1$, then this implies that $X$ and $Y$ are close to each other and we cannot regard $X$ and $Y$ as being *well-separated*. Hence we assume that $\tau$ is a given small or moderate constant (for example, $\tau \leq 0.7$) in the rest of this paper.

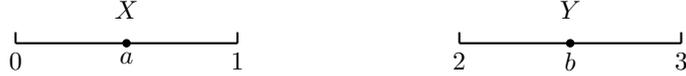

Fig. 2.1: Well-separated intervals $X, Y$ (centered at $a = 0.5, b = 2.5$) with separation ratio $\tau = 0.5$.

Consider the kernel function $1/(x - y)$ again. When $X$ and $Y$ are well-separated so that (2.3) holds, a degenerate expansion for the kernel function based on Taylor expansion takes the following form [18]:

(2.4) $$\kappa(x,y) = \sum_{k=0}^{r-1} \sum_{l=0}^{k} c_{k,l} \phi_{a,l}(x-a) \phi_{b,k-l}(y-b) + \epsilon_r, \quad \forall x \in X, \ y \in Y,$$

where

(2.5)
$$c_{k,l} := \begin{cases} -k!(b-a)^{-(k+1)} \eta_{a,l}^{-1} \eta_{b,k-l}^{-1}(-1)^{k-l} & \text{if } l \leq k, \\ 0 & \text{if } l > k, \end{cases}$$

$$\phi_{v,l}(x) := \eta_{v,l} \frac{x^l}{l!}, \quad \eta_{v,l} = \begin{cases} 1, & \text{if } l = 0, \\ \left(\frac{l}{e}(2\pi r)^{\frac{1}{2r}} \frac{1}{\delta_v}\right)^l & \text{if } l = 1, \ldots, r-1, \end{cases}$$

and the approximation error $\epsilon_r$ satisfies

$$|\epsilon_r| \leq \frac{(1+\tau)\tau^r}{(1-\tau)}|\kappa(x,y)|, \quad \forall x \in X, \ y \in Y. \tag{2.6}$$

The above expansion will be used to illustrate the construction of hierarchical rank structured matrices and analyze the approximation error in the remaining sections.

**2.2. Farfield and nearfield blocks.** We now consider a dense matrix $A$ defined by $A := [\kappa(x,y)]_{x \in X, y \in Y}$. The degenerate approximation (2.4) immediately indicates that certain blocks of $A$ admit a low-rank approximation. In order to identify these low-rank blocks, it is necessary to exploit *nearfield* and *farfield* matrix blocks as they are defined in [12].

DEFINITION 2.2. *Given two sets of points $X_\mathbf{i}$ and $Y_\mathbf{j}$, a submatrix $A|_{\mathbf{i} \times \mathbf{j}}$ is called a* farfield *block if $X_\mathbf{i}$ and $Y_\mathbf{j}$ are well-separated in the sense of Definition 2.1; otherwise, $A|_{\mathbf{i} \times \mathbf{j}}$ is called a* nearfield *block.*

A major difference between farfield and nearfield blocks is that each farfield block can be approximated by low-rank matrices, as a consequence of (2.4). The following theorem restates (2.4) in matrix form.

THEOREM 2.3. *If $X_\mathbf{i}$ and $Y_\mathbf{j}$ are well-separated in the sense of (2.3) with centers $a_i$ and $a_j$, radii $\delta_i$ and $\delta_j$, respectively, the farfield block $A|_{\mathbf{i} \times \mathbf{j}}$ admits a low-rank approximation of the form*

$$A|_{\mathbf{i} \times \mathbf{j}} = \hat{U}_i \hat{B}_{i,j} \hat{V}_j^T + E_F|_{\mathbf{i} \times \mathbf{j}}, \tag{2.7}$$

*where*

$$\hat{U}_i = [\phi_{a_i,l}(x-a_i)]_{\substack{x \in X_\mathbf{i}, \\ l=0:r-1}}, \ \hat{V}_j = [\phi_{a_j,l}(y-a_j)]_{\substack{y \in Y_\mathbf{j}, \\ l=0:r-1}}, \ \hat{B}_{i,j} = [c_{k,l}]_{k,l=0:r-1}, \tag{2.8}$$

*with $c_{k,l}, \phi_{v,l}(v=a_i,a_j)$ defined in (2.5), and*

$$\|E_F|_{\mathbf{i} \times \mathbf{j}}\|_{\max} \leq \epsilon_{far} \|A|_{\mathbf{i} \times \mathbf{j}}\|_{\max} \quad \text{with} \quad \epsilon_{far} = \frac{(1+\tau)\tau^r}{(1-\tau)}. \tag{2.9}$$

Let $n_i = |X_\mathbf{i}|$ and $n_j = |Y_\mathbf{j}|$. If the points $x$ of $X_\mathbf{i}$ are listed as columns and the various functions $\phi_{a_i,l}(x-a_i)$ are listed row-wise with $l = 0, \cdots, r-1$ and similarly for $y$, $Y_\mathbf{j}$, and $\phi_{a_j,l}(y-a_j)$ then the matrices $\hat{U}_i, \hat{B}_{ij}$ and $\hat{V}_j$ has dimensions $n_i \times r$, $r \times r$, and $n_j \times r$, respectively. The theorem is illustrated in Fig.2.2.

**2.3. Strong rank-revealing QR.** Notice that in the approximation (2.7), $\hat{U}_i$ only depends on the points in $X_\mathbf{i}$, $\hat{V}_j$ only depends the points in $Y_\mathbf{j}$ and $\hat{B}_{i,j}$ depends on both the centers of $X_\mathbf{i}$ and $Y_\mathbf{j}$ as well as their radii. This represents a standard expansion structure used in FMM [27, 50, 52]. As will be seen in the next section the construction of $\mathcal{H}^2$ and HSS matrices will be significantly simplified by further postprocessing $\hat{U}_i$ and $\hat{V}_j$ with a strong rank-revealing QR (SRRQR) factorization [30]. The following theorem summarizes Algorithm 4 in [30].

THEOREM 2.4. *([30, Algorithm 4]) Let $M$ be an $m \times n$ matrix and $M \neq 0$. Given a real number $s > 1$ and a positive integer $r$ ($r \leq \text{rank}(A)$), the strong rank-revealing QR algorithm computes a factorization of $M$ in the form:*

$$MP = Q \begin{bmatrix} R_{11} & R_{12} \\ & R_{22} \end{bmatrix}, \tag{2.10}$$

*where $P$ is a permutation matrix, $Q \in \mathbb{R}^{m \times m}$ is an orthogonal matrix, $R_{11}$ is a $r \times r$ upper triangular matrix and $R_{12}$ is a $r \times (n-r)$ dense matrix that satisfies the condition:*

$$\|R_{11}^{-1} R_{12}\|_{\max} \leq s.$$

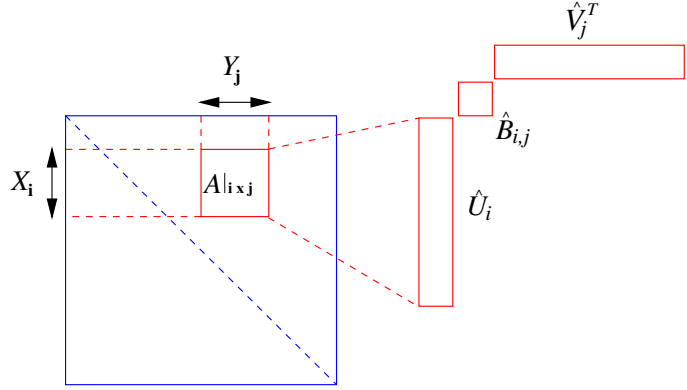

Fig. 2.2: Illustration of Theorem 2.3.

The complexity is $O(n^3 \log_s n)$ if $m \approx n$.

In all of our implementations, we set $s = 2$. SRRQR unravels a set of columns of $A$ that nearly span the range of $A$ – thus the term rank-revealing. Assume $C$ is an $n \times r$ matrix with rank $r$. Applying SRRQR to $C^T$ produces the following factorization:

$$C^T P = Q \begin{bmatrix} R_{11} & R_{12} \end{bmatrix},$$

where $Q \in \mathbb{R}^{r \times r}$ is an orthogonal matrix. A modification of the above equation leads to

$$C = P \begin{bmatrix} I \\ (R_{11}^{-1} R_{12})^T \end{bmatrix} \tilde{C},$$

where $I$ is an identity matrix of order $r$ and $\tilde{C} = (QR_{11})^T$. Note that the above relation implies that $\tilde{C} \in \mathbb{R}^{r \times r}$ is a submatrix of $C$ consisting of the first $r$ rows of the row-permuted matrix $P^T C$. From this perspective, the whole aim of the procedure is to extract a set of $r$ rows from $C$ that will nearly span its row space.

When $\hat{U}_i$ and $\hat{V}_j$ in (2.7) both have more rows than columns, applying SRRQR to $\hat{U}_i^T$ and then to $\hat{V}_j^T$ yields:

$$(2.11) \qquad \hat{U}_i = P_i \begin{bmatrix} I \\ G_i \end{bmatrix} \hat{U}_i|_{\hat{\mathbf{i}}}, \quad \hat{V}_j = F_j \begin{bmatrix} I \\ H_j \end{bmatrix} \hat{V}_j|_{\hat{\mathbf{j}}}.$$

Note that, as explained above for $\tilde{C}$, $\hat{U}_i|_{\hat{\mathbf{i}}}$ denotes a matrix made up of selected rows of $\hat{U}_i$.

Substituting the above two equations into (2.7) leads to another form of the low-rank approximation to $A|_{\mathbf{i} \times \mathbf{j}}$:

$$(2.12) \qquad A|_{\mathbf{i} \times \mathbf{j}} \approx P_i \begin{bmatrix} I \\ G_i \end{bmatrix} \hat{U}_i|_{\hat{\mathbf{i}}} \hat{B}_{i,j} (\hat{V}_j|_{\hat{\mathbf{j}}})^T \left( F_j \begin{bmatrix} I \\ H_j \end{bmatrix} \right)^T$$

$$(2.13) \qquad = P_i \begin{bmatrix} I \\ G_i \end{bmatrix} (A|_{\hat{\mathbf{i}} \times \hat{\mathbf{j}}} - E_F|_{\hat{\mathbf{i}} \times \hat{\mathbf{j}}}) \left( F_j \begin{bmatrix} I \\ H_j \end{bmatrix} \right)^T \approx P_i \begin{bmatrix} I \\ G_i \end{bmatrix} A|_{\hat{\mathbf{i}} \times \hat{\mathbf{j}}} \left( F_j \begin{bmatrix} I \\ H_j \end{bmatrix} \right)^T,$$

where $\hat{\mathbf{i}}$ and $\hat{\mathbf{j}}$ represent subsets of $\mathbf{i}$ and $\mathbf{j}$, respectively, and (2.13) results from (2.7).

A major advantage of this form of approximation over (2.7) is a reduction in storage. Now only four index sets are needed to represent $(P_i, F_j, \hat{\mathbf{i}}, \hat{\mathbf{j}})$ and two smaller dense matrices $(G_i, H_j)$ need to be stored rather than three dense matrices ($\hat{U}_i, \hat{B}_{i,j}, \hat{V}_j$). There are other advantages that will be discussed in the next section.

Fig. 3.1: Illustration of an adaptive partitioning for the case $X = Y = \{x_1, x_2, \ldots, x_8\}$. Left: the computational domain $\Omega$ is recursively bisected until the number of points in each sub-interval $\Omega_i$ centered at $a_i$ is less than 4 (circled dots represent the points $x_i$). Right: the corresponding postordered binary tree $\mathcal{T}$ with indices of points stored at each node.

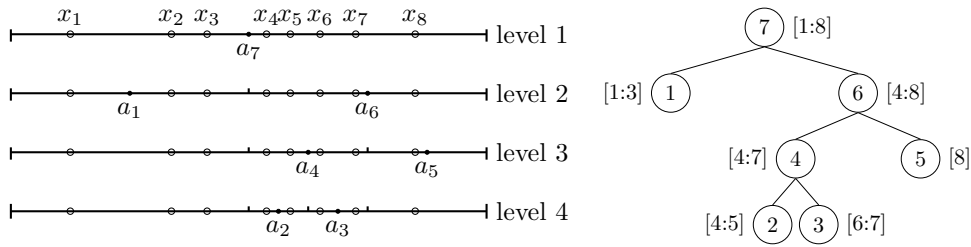

The operations represented by (2.11) will be used extensively in the construction of hierarchical matrices to be seen in the next section. These will be denoted as follows:

$$(2.14) \qquad [P_i, G_i, \hat{\mathbf{i}}] = \text{compr}(\hat{U}_i, \mathbf{i}) \quad \text{and} \quad [F_j, H_j, \hat{\mathbf{j}}] = \text{compr}(\hat{V}_j, \mathbf{j}).$$

Each of the above operations is also called a structure-preserving rank-revealing (SPRR) factorization [60] or an interpolative decomposition [36]. For recent developments on rank-revealing algorithms, see [28].

**3. Construction of hierarchical rank structured matrices with nested bases.** This section presents SMASH, an algorithm to construct either an $\mathcal{H}^2$ or an HSS matrix approximation to an $n \times n$ matrix $A := [\kappa(x, y)]_{x \in X, y \in Y}$, where $\kappa$ is a given kernel function and $X$ and $Y$ are two finite sets of points. Although the discussion focuses on square matrices, SMASH can be extended to rectangular ones [54] without any difficulty.

**3.1. Adaptive partitioning.** SMASH starts with a hierarchical partitioning of the computational domain $\Omega$ and then builds a tree structure $\mathcal{T}$ to facilitate subsequent operations. In order to deal with the case when $X$ and $Y$ are non-uniformly distributed, an adaptive partitioning scheme is necessary.

Without loss of generality, assume both $X$ and $Y$ are contained in a unit box $\Omega = [0, 1]^d$ in $\mathbb{R}^d$ ($d = 1, 2, 3$). The basic idea of this partitioning algorithm is to recursively subdivide the computational domain $\Omega$ into several subdomains until the number of points included in each resulting subdomain is less than a prescribed constant $\nu_0$ (usually much smaller than the number of points in the domain). Specifically, at level 1, $\Omega$ is not partitioned. Starting from level $l$ ($l \geq 2$), each box obtained at level $l-1$ that contains more than $\nu_0$ points is bisected along each of the $d$ dimensions.

For convenience we assume that the number of points from $X$ and $Y$ in each partitioning is the same. If a box is empty, it is discarded. Let $L$ be the maximum level where the recursion stops. Then the information about the partitioning can be represented by a tree $\mathcal{T}$ with $L$ levels, where the root node is at level 1 and corresponds to the domain $\Omega$ and each nonroot node $i$ corresponds to a partitioned subdomain $\Omega_i$. See Fig. 3.1 for a 1D example. The adaptive partitioning guarantees that each subdomain corresponding to a leaf node contains a small number of points less than the prescribed constant $\nu_0$.

**3.2. Review of $\mathcal{H}^2$ and HSS matrices.** The low-rank property of a block $A|_{\mathbf{i}\times\mathbf{j}}$ associated with a node pair $(i,j)$ is related to the *strong (or standard) admissibility condition* employed to define $\mathcal{H}$ and $\mathcal{H}^2$ matrices ([14],[35]):

> for a fixed $\tau \in (0,1)$, the node pair $(i,j)$ in $\mathcal{T}$ is admissible if $X_\mathbf{i}$ and $Y_\mathbf{j}$ are well-separated in the sense of Definition 2.1.

Hierarchical matrices are often defined in terms of the above condition, which, in essence, spells out when a given block in the matrix can be compressed. A matrix $\hat{A}$ (associated with a tree $\mathcal{T}$) is called an $\mathcal{H}$ matrix ([33]) of rank $r$ if there exists a positive integer $r$ such that

$$\text{rank}(\hat{A}|_{\mathbf{i}\times\mathbf{j}}) \leq r, \quad \text{whenver } (i,j) \text{ is admissible.}$$

Furthermore, $\hat{A}$ is called a *uniform $\mathcal{H}$ matrix* ([33]) if there exist a column basis set $\{U_i\}_{i\in\mathcal{T}}$ and a row basis set $\{V_i\}_{i\in\mathcal{T}}$ associated with $\mathcal{T}$, such that when $(i,j)$ is admissible, $\hat{A}|_{\mathbf{i}\times\mathbf{j}}$ admits a low-rank factorization:

$$\hat{A}|_{\mathbf{i}\times\mathbf{j}} = U_i B_{i,j} V_j^T, \quad \text{for some matrix } B_{i,j}.$$

This factorization is referred to as an *AKB representation* in [31], where $B_{i,j}$ is called a *coefficient matrix*. In [12], $B_{i,j}$ is termed a *coupling matrix* and we will follow this terminology here.

The class of $\mathcal{H}^2$ matrices [33] is a subset of uniform $\mathcal{H}$ matrices with a more refined structure. That is, $\hat{A}$ is an $\mathcal{H}^2$ matrix if it is a uniform $\mathcal{H}$ matrix with nested bases in the sense that one can readily express a basis at one level from that of its children (see (3.2)). What is exploited here is that admissible blocks are low-rank and in addition their factors (or generators) can be expressed from lower levels.

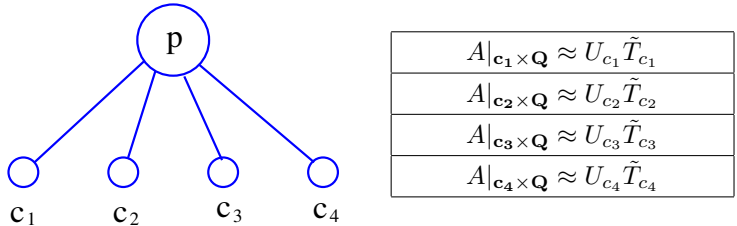

Fig. 3.2: A parent node $p$ with children $c_1, \cdots, c_4$ and the corresponding partition of the matrix block $A|_{\mathbf{p}\times\mathbf{Q}}$, where $\mathbf{Q}$ is the collection of indices associated with all nodes $q$ such that $(p,q)$ is admissible. In the context of HSS matrices there are at most 2 children since the trees are binary. For $\mathcal{H}^2$ matrices the trees are more general.

Assume we have a situation illustrated in Fig. 3.2 where the parent node $p$ has children $c_1, \ldots, c_4$. According to the interpolation in (2.2), the column basis $U_i$ associated with the set $\mathbf{i}$ (for any nonroot node $i$) can be chosen as

$$U_i = \left[p_k^{(i)}(x)\right]_{\substack{x\in X_\mathbf{i} \\ k=1:r}},$$

where $p_k^{(i)}$ ($k = 1, \ldots, r$) are Lagrange basis polynomials corresponding to interpolation points $x_{i_1}, \ldots, x_{i_r}$. If $i$ is a child of node $p$, we can write (cf.[33])

(3.1) $$p_k^{(p)}(x) = \sum_{l=1:r} p_k^{(p)}(x_{i_l}) p_l^{(i)}(x).$$

The matrix version of (3.1) then leads to the so-called *nested basis*:

$$U_p = \begin{bmatrix} U_{c_1} R_{c_1} \\ \vdots \\ U_{c_4} R_{c_4} \end{bmatrix}, \quad \text{with} \quad R_i = \begin{bmatrix} p_1^{(p)}(x_{i_1}) & \cdots & p_r^{(p)}(x_{i_1}) \\ \vdots & \vdots & \vdots \\ p_1^{(p)}(x_{i_r}) & \cdots & p_r^{(p)}(x_{i_r}) \end{bmatrix}.$$

The nested basis can also be obtained through algebraic compressions based on a bottom-up procedure. Let $A|_{\mathbf{p} \times \mathbf{Q}}$ denote the entire (numerically) low rank block row associated with node $p$, i.e., $\mathbf{Q}$ is the union of all indices $\mathbf{q}$ such that $(p,q)$ is admissible. As illustrated in Fig. 3.2, assuming that the column basis $U_{c_i}$ has been obtained from a rank-$r$ factorization of the submatrix $A|_{\mathbf{c_i} \times \mathbf{Q}} \approx U_{c_i} \tilde{T}_{c_i}$, we then derive

$$A|_{\mathbf{p} \times \mathbf{Q}} \approx \begin{bmatrix} U_{c_1} & & & \\ & U_{c_2} & & \\ & & U_{c_3} & \\ & & & U_{c_4} \end{bmatrix} \begin{bmatrix} \tilde{T}_{c_1} \\ \tilde{T}_{c_2} \\ \tilde{T}_{c_3} \\ \tilde{T}_{c_4} \end{bmatrix}.$$

Applying a rank-$r$ factorization to the transpose of $\begin{bmatrix} \tilde{T}_{c_1}^T & \tilde{T}_{c_2}^T & \tilde{T}_{c_3}^T & \tilde{T}_{c_4}^T \end{bmatrix}$ leads to

$$\begin{bmatrix} \tilde{T}_{c_1} \\ \tilde{T}_{c_2} \\ \tilde{T}_{c_3} \\ \tilde{T}_{c_4} \end{bmatrix} \approx \begin{bmatrix} R_{c_1} \\ R_{c_2} \\ R_{c_3} \\ R_{c_4} \end{bmatrix} T_p \longrightarrow A|_{\mathbf{p} \times \mathbf{Q}} \approx U_p T_p \quad \text{with} \quad U_p = \begin{bmatrix} U_{c_1} R_{c_1} \\ U_{c_2} R_{c_2} \\ U_{c_3} R_{c_3} \\ U_{c_4} R_{c_4} \end{bmatrix}.$$

Thus, we can get the basis $U_p$ for the parent node from the children's $U_{c_i}$'s and the matrices $R_{c_i}$ from both analytic and algebraic compression schemes. The $R_i$s are called transfer matrices. Clearly, a similar process can be applied to obtain a row-basis $V_p$ and so, more generally, we can write

$$(3.2) \qquad U_p = \begin{bmatrix} U_{c_1} R_{c_1} \\ \vdots \\ U_{c_k} R_{c_k} \end{bmatrix}, \quad V_p = \begin{bmatrix} V_{c_1} W_{c_1} \\ \vdots \\ V_{c_k} W_{c_k} \end{bmatrix}.$$

Hence only the matrices $U_i, V_i$ for all leaf nodes $i$ need to be stored. Matrices $U_p, V_p$ for a non-leaf node $p$ can be obtained via transfer matrices which require much less storage. This is at the origin of the improvement from an $O(n \log n)$ cost for the early method in this class developed by Barnes and Hut [3] ($\mathcal{H}$ structure) into an $O(n)$ cost method obtained later by the FMM [27] ($\mathcal{H}^2$ structure) for computing matrix-vector multiplications for some kernel matrices ([52]).

Note that as they are described in the literature $\mathcal{H}$ and $\mathcal{H}^2$ matrices are associated with more general trees than those traditionally used for HSS matrices [21, 22] which are binary trees, according to the partitioning algorithm described in Section 3.1. In fact HSS matrices can be viewed as a special class of $\mathcal{H}^2$ matrices in which the strong admissibility condition is replaced by the *weak admissibility condition*[35]:

the node pair $(i, j)$ in $\mathcal{T}$ is admissible if $i \neq j$.

The above weak admissibility condition implies that, if $\hat{A}$ is an HSS matrix, and $i, j$ are two children of the root node, then the matrix block $\hat{A}|_{\mathbf{i} \times \mathbf{j}}$ should admit a low-rank factorization.

In the context of integral equations, this requirement means that the HSS structure will face difficulties in situations when couplings between nearfield blocks require a relatively high rank representation. Approximation by HSS matrices will work well for integral equations defined on a curve where kernel functions are discretized. In other cases the numerical rank of $A|_{\mathbf{i}\times\mathbf{j}}$ may not necessarily be small even when a non-oscillatory kernel function is discretized on a surface or in a volume [12, 56].

The construction of $\mathcal{H}^2$ and HSS matrices involves computing the basis matrices $U, V$ at the leaf level, along with the transfer matrices $R, W$, and the coupling matrices $B$ associated with a tree $\mathcal{T}$. In particular, each leaf node $i$ is assigned four matrices $\{U_i, V_i, R_i, W_i\}$ and each nonleaf node $i$ at level $l \geq 3$ is assigned two matrices $\{R_i, W_i\}$.

There are two types of $B_{i,j}$ matrices, those corresponding to the nearfield blocks at the leaf level, and those corresponding to the coupling matrix where the product $U_i B_{i,j} V_j^T$ approximates block $A|_{\mathbf{i}\times\mathbf{j}}$ for *certain* admissible $(i,j)$. In general, the computation of $B_{i,j}$ is more complicated because one has to carefully specify the set of admissible node pairs $(i,j)$ to be used for the efficient approximation of $A$. If the distribution of points is uniform, the corresponding node pairs $(i,j)$ are related to what is called *interaction list* in FMM [27, 52]. In more general settings where points can be arbitrarily distributed, they are called *admissible leaves* [12]. The set of admissible leaves corresponding to the minimal admissible partitioning [31] can be defined as follows:

(3.3)
$$\begin{aligned}
\mathcal{L} = & \{(i,j) : i,j \in \mathcal{T} \text{ are nodes at the same level such that } (i,j) \text{ is admissible} \\
& \text{but } (p_i, p_j) \text{ is not admissible}, \text{ where } p_i, p_j \text{ are parents of } i, j, \text{ respectively}\} \\
& \cup \{(i,j) : i \in \mathcal{T} \text{ is a leaf node and } j \in \mathcal{T} \text{ with } \text{lv}(j) > \text{lv}(i) \text{ such that} \\
& (i,j) \text{ is admissible but } (i, p_j) \text{ is not admissible with } p_j \text{ the parent of } j\} \\
& \cup \{(i,j) : j \in \mathcal{T} \text{ is a leaf node and } i \in \mathcal{T} \text{ with } \text{lv}(i) > \text{lv}(j) \text{ such that} \\
& (i,j) \text{ is admissible but } (p_i, j) \text{ is not admissible with } p_i \text{ the parent of } i\}.
\end{aligned}$$

The node pairs $(i,j)$ corresponding to blocks $B_{i,j}$ that can not be compressed or partitioned, can be identified through inadmissible leaves as defined below (cf.[12]):

(3.4) $$\mathcal{L}^- := \{(i,j) : i,j \in \mathcal{T} \text{ are leaf nodes and } (i,j) \text{ is not admissible}\}.$$

In particular, for HSS matrices, it can be seen that $\mathcal{L}$ and $\mathcal{L}^-$ have the following simple expression:

(3.5) $$\mathcal{L} = \{(i,j) : i,j \in \mathcal{T} \text{ and } j = \text{sibling of } i\}, \quad \mathcal{L}^- = \{(i,i) : i \in \mathcal{T} \text{ is a leaf node}\}.$$

This special feature will be used later (Section 3.3.1) to simplify the notation associated with HSS matrices. The $U, V, R, W, B$ matrices are called $\mathcal{H}^2$, or HSS, *generators* in the remaining sections.

**3.3. Levelwise construction.** Although the HSS structure may appear to be simpler than the $\mathcal{H}^2$ structure, based on their algebraic definitions the HSS construction procedure is actually more complex. This is because HSS matrices require the compression of both nearfield and farfield blocks while $\mathcal{H}^2$ matrices only require the compression on farfield blocks. For example, if two sets $X_{\mathbf{i}}$ and $Y_{\mathbf{j}}$ are almost adjacent to each other ($\tau \approx 1$ in (2.3)), then the analytic approximation will not produce a *low rank*, i.e., to get an accurate approximation, $r$ has to be large in (2.6). In this case, the $\mathcal{H}^2$ matrix will form this block explicitly as a dense matrix while the HSS matrix still requires the block to be factorized. In what follows, we will first discuss SMASH for the HSS construction in detail and then present the $\mathcal{H}^2$ construction with an emphasis on their differences.

**3.3.1. HSS construction.** Due to the simple structure of $\mathcal{L}, \mathcal{L}^-$ in (3.5), the notation denoting the coupling matrices and nearfield blocks can be simplified in the HSS representation. Specifically, for $(i,j) \in \mathcal{L}$, $B_{i,j}$ can be represented as $B_i$ with the second index $j$ dropped because $j = $ sibling of $i$ is uniquely determined in a binary tree. An additional symbol $D_i$ is introduced to represent diagonal blocks $B_{i,i}$ because $(i,j) \in \mathcal{L}^-$ implies $j = i$.

The basic idea of SMASH for the HSS construction is to first apply a truncated SVD to obtain a basis for nearfield blocks, use interpolation or expansion to obtain a basis for farfield blocks and then apply SRRQR to the combination of these two bases to obtain the $U, V, R$ and $W$ matrices. The $D$ and $B$ matrices are submatrices of the original kernel matrix and their indices are readily available after the computation of $U, V, R, W$ matrices. In order to distinguish between column and row indices associated with a node $i$, we use superscript $row$ to mark its row indices and $col$ to mark its column indices. For example, $\mathbf{i}^{row}$ and $\mathbf{i}^{col}$ denote the indices of points from X and Y contained in $\Omega_i$, respectively.

Assuming the HSS tree $\mathcal{T}$ has $L$ levels, the HSS construction algorithm traverses $\mathcal{T}$ through level $l = L, L-1, \ldots, 2$. Before the construction, two intermediate sets $\bar{\mathbf{i}}^{row}$ and $\bar{\mathbf{i}}^{col}$ are initialized as follows for each node $i$:

$$\bar{\mathbf{i}}^{row} = \begin{cases} \mathbf{i}^{row} & \text{if } i \text{ is a leaf,} \\ \emptyset & \text{otherwise,} \end{cases} \quad \bar{\mathbf{i}}^{col} = \begin{cases} \mathbf{i}^{col} & \text{if } i \text{ is a leaf,} \\ \emptyset & \text{otherwise,} \end{cases} \tag{3.6}$$

where the index sets $\mathbf{i}^{row}$ and $\mathbf{i}^{col}$ have been saved for each leaf node after the partitioning of $\Omega$.

Let $\mathcal{N}_i$ be the set of blocks that are nearfield to node $i$. This set is used for HSS matrices only and is defined in the appendix. For each node $i$ at level $l$, the construction algorithm first applies a truncated SVD to compute an approximate column basis for the nearfield block rows in terms of $X_{\bar{\mathbf{i}}^{row}}$:

$$A_i^- := \left[A|_{\bar{\mathbf{i}}^{row} \times \bar{\mathbf{j}}^{col}}\right]_{j \in \mathcal{N}_i} = S_i \Sigma_i^- \left[\widetilde{S}_j\right]_{j \in \mathcal{N}_i} + \left[E_\Sigma^-|_{\bar{\mathbf{i}}^{row} \times \bar{\mathbf{j}}^{col}}\right]_{j \in \mathcal{N}_i}, \tag{3.7}$$

where the columns of $S_i$ and $[\widetilde{S}_j]_{j \in \mathcal{N}_i}$ are the left/right singular vectors of $A_i^-$ and $\Sigma_i^-$ is a diagonal matrix composed of corresponding singular values of $A_i^-$ such that the following estimate holds

$$\|E_\Sigma^-|_{\bar{\mathbf{i}}^{row} \times \bar{\mathbf{j}}^{col}}\|_F \leq \sqrt{|\bar{\mathbf{i}}^{row}||\bar{\mathbf{j}}^{col}|} \epsilon_{\text{SVD}} \|A_i^-\|_2, \quad \forall j \in \mathcal{N}_i. \tag{3.8}$$

Here, $\epsilon_{\text{SVD}}$ is the relative tolerance used in the truncated SVD. This estimate will be used to analyze the overall construction error in Section 4. The matrix $S_i$ is then taken as an approximate column basis for the nearfield block rows $A_i^-$.

For farfield blocks with respect to $X_{\bar{\mathbf{i}}^{row}}$, a column basis $\hat{U}_i$ can be easily obtained through interpolation (2.2) or Taylor expansion (2.8) that only rely on $X_{\bar{\mathbf{i}}^{row}}$ and $\Omega_i$. Next, we apply SRRQR to the combined basis $[\hat{U}_i, S_i]$ as shown below

$$[P_i, G_i, \hat{\mathbf{i}}^{row}] = \text{compr}([\hat{U}_i, S_i], \bar{\mathbf{i}}^{row}). \tag{3.9}$$

From these outputs, we set

$$\begin{aligned} U_i &:= P_i \begin{bmatrix} I \\ G_i \end{bmatrix} \quad \text{if } i \text{ is a leaf node,} \\ \begin{bmatrix} R_{c_1} \\ R_{c_2} \end{bmatrix} &:= P_i \begin{bmatrix} I \\ G_i \end{bmatrix} \quad \text{if } i \text{ is a parent with children } c_1, c_2. \end{aligned} \tag{3.10}$$

Similarly, in order to compute $V, W$ generators, a truncated SVD is first applied to the nearfield block columns (transposed) in terms of $Y_{\bar{\mathbf{j}}^{col}}$:

$$A_i^| := \left[\left(A_{\bar{\mathbf{j}}^{row} \times \bar{\mathbf{i}}^{col}}\right)^T\right]_{j \in \mathcal{N}_i} = T_i \Sigma_i^| \left[\widetilde{T}_j\right]_{j \in \mathcal{N}_i} + \left[\left(E_\Sigma^||_{\bar{\mathbf{j}}^{row} \times \bar{\mathbf{i}}^{col}}\right)^T\right]_{j \in \mathcal{N}_i}, \tag{3.11}$$

where the truncation error satisfies

$$(3.12) \qquad \|E_\Sigma^|\,|_{\bar{\mathbf{j}}^{row}\times\bar{\mathbf{i}}^{col}}\|_F \leq \sqrt{|\bar{\mathbf{i}}^{col}||\bar{\mathbf{j}}^{row}|}\epsilon_{\text{SVD}}\|A_i^|\|_2, \quad \forall\, j \in \mathcal{N}_i.$$

In the nex step we compute a row basis $\hat{V}_i$ for the farfield blocks with respect to $Y_{\bar{\mathbf{i}}^{col}}$ based on (2.2) or (2.8) and apply SRRQR to $[\hat{V}_i, T_i]$:

$$[F_i, H_i, \hat{\mathbf{i}}^{col}] = \text{compr}([\hat{V}_i, T_i], \bar{\mathbf{i}}^{col}).$$

Then we set

$$(3.13) \qquad \begin{aligned} V_i &:= F_i \begin{bmatrix} I \\ H_i \end{bmatrix} \quad \text{if } i \text{ is a leaf node,} \\ \begin{bmatrix} W_{c_1} \\ W_{c_2} \end{bmatrix} &:= F_i \begin{bmatrix} I \\ H_i \end{bmatrix} \quad \text{if } i \text{ is a parent with children } c_1, c_2. \end{aligned}$$

Once the compressions for children nodes (at level $l$) are complete, we update the intermediate index set associated with the parent node (at level $l-1$) as shown below :

$$(3.14) \qquad \bar{\mathbf{p}}^{row} = \hat{\mathbf{c}_1}^{row} \cup \hat{\mathbf{c}_2}^{row}, \quad \bar{\mathbf{p}}^{col} = \hat{\mathbf{c}_1}^{col} \cup \hat{\mathbf{c}_2}^{col},$$

where $c_1, c_2$ are the children of $p$.

After the bottom-up traversal of $\mathcal{T}$ and hence the computation of $U, V, R, W$ matrices, the $B$ and $D$ matrices can be extracted as follows:

$$(3.15) \qquad B_i := A|_{\hat{\mathbf{i}}^{row}\times\hat{\mathbf{j}}^{col}}, \quad j = \text{sibling of } i, \quad \text{and} \quad D_i := A|_{\mathbf{i}^{row}\times\mathbf{i}^{col}}, \quad i = \text{leaf node}.$$

**3.3.2. $\mathcal{H}^2$ construction.** As mentioned at the beginning of Section 3.3, the $\mathcal{H}^2$ construction is simpler because nearfield blocks will not be factorized, and the only complication is that an $\mathcal{H}^2$ matrix may be associated with a more general tree structure where a parent can have more than two children.

SMASH for the $\mathcal{H}^2$ construction also follows a bottom-up levelwise traversal of $\mathcal{T}$ through level $l = L, L-1, \ldots, 2$. For each node $i$ at level $l$, a column/row basis $\hat{U}_i/\hat{V}_i$ corresponding to a farfield block row/column with index $\bar{\mathbf{i}}^{row}/\bar{\mathbf{i}}^{col}$ can be obtained via either interpolation (2.2) or Taylor expansion (2.8). The bases $\hat{U}_i$ and $\hat{V}_i$ are then passed into SRRQR

$$(3.16) \qquad [P_i, G_i, \hat{\mathbf{i}}^{row}] = \text{compr}(\hat{U}_i, \bar{\mathbf{i}}^{row}) \quad \text{and} \quad [F_i, H_i, \hat{\mathbf{i}}^{col}] = \text{compr}(\hat{V}_i, \bar{\mathbf{i}}^{col}).$$

The $\mathcal{H}^2$ generators $U, R, V, W$ are then set as

$$(3.17) \qquad \begin{aligned} U_i &:= P_i \begin{bmatrix} I \\ G_i \end{bmatrix} \quad \text{if } i \text{ is a leaf node,} \\ \begin{bmatrix} R_{c_1} \\ \vdots \\ R_{c_k} \end{bmatrix} &:= P_i \begin{bmatrix} I \\ G_i \end{bmatrix} \quad \text{if } i \text{ is a parent with children } c_1, \ldots, c_k, \\ V_i &:= F_i \begin{bmatrix} I \\ H_i \end{bmatrix} \quad \text{if } i \text{ is a leaf node,} \\ \begin{bmatrix} W_{c_1} \\ \vdots \\ W_{c_k} \end{bmatrix} &:= F_i \begin{bmatrix} I \\ H_i \end{bmatrix} \quad \text{if } i \text{ is a parent with children } c_1, \ldots, c_k. \end{aligned}$$

Again, once the compressions for children nodes (at level $l$) are complete, the intermediate index set associated with the parent node (at level $l-1$) can be updated as in (3.14). Namely,

$$(3.18) \quad \bar{\mathbf{p}}^{row} = \hat{\mathbf{c_1}}^{row} \cup \cdots \cup \hat{\mathbf{c_k}}^{row}, \quad \bar{\mathbf{p}}^{col} = \hat{\mathbf{c_1}}^{col} \cup \cdots \cup \hat{\mathbf{c_k}}^{col}, \quad \text{with children}(p) = \{c_1, \ldots, c_k\}.$$

Finally, analogous to (3.15), the coupling matrices associated with admissible leaves are extracted based on index sets $\hat{\mathbf{i}}^{row}$ and $\hat{\mathbf{j}}^{col}$ as

$$(3.19) \quad B_{i,j} := A|_{\hat{\mathbf{i}}^{row} \times \hat{\mathbf{j}}^{col}}, \quad \forall (i,j) \in \mathcal{L},$$

and the nearfield blocks associated with inadmissible leaves are formed by

$$(3.20) \quad B_{i,j} := A|_{\mathbf{i}^{row} \times \mathbf{j}^{col}}, \quad \forall (i,j) \in \mathcal{L}^-.$$

Compared with standard $\mathcal{H}^2$ constructions based on either expansion or interpolation, SMASH is more efficient and easier to implement. First, in order to complete the $\mathcal{H}^2$ construction procedure, it suffices to provide only the column/row basis for each farfield block, which can be easily obtained based on interpolation (2.2), for example, and the coupling matrices $B_{i,j}$ can be simply extracted from the original matrix without resorting to any other formulas. Second, no information is required about the translation to compute transfer matrices because the computation of $R/W$ is essentially the same as that of $U/V$ at leaf level. For all the children $i$ of a node $p$, $R_i/W_i$ are calculated jointly based on a subset of points located inside $\Omega_p$ (i.e., $X_{\bar{\mathbf{p}}^{row}}/Y_{\bar{\mathbf{p}}^{col}}$). Therefore, SMASH essentially builds a hierarchy of grids and computes the bases at each level of the tree by repeating the same operations (3.16) on each coarse grid (this can be more clearly seen in (4.1) in view of a perfect tree). In addition, the use of SRRQR guarantees that each entry of the $U, V, R, W$ matrices is bounded by a user-specified constant, which ensures the numerical stability of the construction procedure. Note that the special structures in the nested bases (3.17) result in not only a reduced storage but also in faster matrix operations such as matrix-vector multiplications, linear system solutions, etc. Finally, since the computation of the basis matrices only relies on the information local to each node, as can be seen from (3.9) and (3.16), this construction algorithm is inherently suitable for a parallel computing environment.

**3.4. Matrix-vector multiplication.** Among various hierarchical rank structured matrix operations, the matrix-vector multiplication is the most widely used, as indicated by the popularity of tree code [3] (for $\mathcal{H}$ matrices) and fast multipole method [27, 52] (for $\mathcal{H}^2$ matrices).

The matrix-vector multiplication for an $\mathcal{H}^2$ matrix $A$ follows first a bottom-up and then a top-down traversal of $\mathcal{T}$ [12, 31], which is a succinct algebraic generalization of the fast multipole method (cf.[52]). Suppose $\mathcal{T}$ has $L$ levels, the node-wise version of this algorithm to evaluate $z = Aq$ can be summarized as follows:

1. from level $l = L$ to level $l = 2$, for each node $i$ at level $l$, compute $\hat{q}_i := V_i^T q|_{\mathbf{i}^{col}}$ if $l = L$; otherwise, compute $\hat{q}_i := \sum_{c \in \text{children}(i)} W_c^T \hat{q}_c$;
2. for each nonroot node $i \in \mathcal{T}$, compute $\hat{z}_i = \sum_{j:(i,j) \in \mathcal{L}} B_{i,j} \hat{q}_j$;
3. from level $l = 2$ to level $l = L$, for each node $i$ at level $l$, if $l < L$, for each child $c$ of $i$, compute $\hat{z}_c = \hat{z}_c + R_c \hat{z}_i$; otherwise, compute $z|_{\mathbf{i}^{row}} = U_i \hat{z}_i + \sum_{j:(i,j) \in \mathcal{L}^-} B_{i,j} q|_{\mathbf{j}^{col}}$.

When $X$ and $Y$ are uniformly distributed in $[0,1]^d$, the resulting tree $\mathcal{T}$ is a perfect $2^d$-tree (each parent has $2^d$ children and all leaf nodes are at the same level). If, in addition we assume the ordering of points to be consistent with the postordering of tree $\mathcal{T}$, i.e., for two siblings $i, j \in \mathcal{T}$, if $i < j$ then the index of any point in box $i$ must be smaller than point in box $j$, then an $\mathcal{H}^2$ matrix $A$ has a telescoping representation:

$$A = B^{(L)} +$$
$$(3.21) \quad U^{(L)} \left( U^{(L-1)} \left( \ldots \left( U^{(2)} B^{(1)} \left( V^{(2)} \right)^T + B^{(2)} \right) \ldots \right) \left( V^{(L-1)} \right)^T + B^{(L-1)} \right) \left( V^{(L)} \right)^T,$$

where $U^{(l)}, V^{(l)}$ are block diagonal matrices:

$$(3.22) \quad U^{(l)} = \begin{cases} \text{diag}(U_i)_{\text{lv}(i)=l} & \text{if } l = L, \\ \text{diag}\left(\begin{bmatrix} R_{c_1} \\ \vdots \\ R_{c_k} \end{bmatrix}\right)_{\text{lv}(i)=l} & \text{if } l < L, \end{cases} \quad V^{(l)} = \begin{cases} \text{diag}(V_i)_{\text{lv}(i)=l} & \text{if } l = L, \\ \text{diag}\left(\begin{bmatrix} W_{c_1} \\ \vdots \\ W_{c_k} \end{bmatrix}\right)_{\text{lv}(i)=l} & \text{if } l < L, \end{cases}$$

and $B^{(l)}$ has a block structure. $B^{(L)}$ has $\#\{i \in \mathcal{T} : \text{lv}(i) = L\} \times \#\{i \in \mathcal{T} : \text{lv}(i) = L\}$ blocks where each nonzero block corresponds to a nearfield block, while for $l < L$, there are $\#\{i \in \mathcal{T} : \text{lv}(i) = l+1\} \times \#\{i \in \mathcal{T} : \text{lv}(i) = l+1\}$ blocks in $B^{(l)}$ and each nonzero block corresponds to a coupling matrix. That is, block $(l_i, l_j)$ is equal to $B_{i,j}$ if $\text{lv}(i) = \text{lv}(j) = l$ such that $(i,j) \in \mathcal{L}^-$ when $l = L$ or $(i,j) \in \mathcal{L}$ when $l \leq L$. Here $l_i$ denotes the location of node $i$ at level $l$ enumerated from left to right. If $A$ is an HSS matrix associated with a perfect binary tree $\mathcal{T}$, the structures of $U^{(l)}$ and $V^{(l)}$ are identical to those in (3.22) with $k = 2$ but $B^{(l)}$ has a much simpler block diagonal structure:

$$(3.23) \quad B^{(l)} = \begin{cases} \text{diag}(D_i)_{i \text{ is a leaf node}} & \text{if } l = L, \\ \text{diag}\left(\begin{bmatrix} 0 & B_{c_1} \\ B_{c_2} & 0 \end{bmatrix}\right)_{\text{lv}(i)=l} & \text{if } l < L, \end{cases}$$

where $c_1$ and $c_2$ are the children of node $i$.

Based on the explicit representation (3.21) of an $\mathcal{H}^2$ matrix associated with a perfect $2^d$-tree, we can write down a levelwise version of the matrix-vector multiplication:

1. at level $l (2 \leq l \leq L)$, compute

$$(3.24) \quad \hat{q}^{(l)} = \left(V^{(l)}\right)^T \ldots \left(V^{(L)}\right)^T q;$$

2. at level $l (2 \leq l \leq L)$, compute

$$(3.25) \quad \hat{z}^{(l)} = B^{(l-1)} \hat{q}^{(l)};$$

3.

$$(3.26) \quad z = B^{(L)} q + U^{(L)} \left( \ldots \left( U^{(3)} (U^{(2)} \hat{z}^{(2)} + \hat{z}^{(3)}) + \hat{z}^{(4)} \right) + \ldots \hat{z}^{(L)} \right).$$

Notice that when $X, Y$ are not uniformly distributed, $\mathcal{T}$ is not necessarily a perfect tree. Under this condition, the nodes $i, j$ corresponding to a coupling matrix $B_{i,j}$ may not be at the same level of $\mathcal{T}$ and the telescoping expansion (3.21) does not exist.

As for linear system solutions, $\mathcal{H}^2$ and HSS matrices take completely different approaches to directly solve the resulting system. Linear complexity $\mathcal{H}^2$ matrix solvers ([12, 31]) heavily depend on recursion to reduce the number of floating point operations while HSS matrices could benefit from highly parallelizable ULV-type algorithms (cf.[21], [22]) due to the special structure of HSS. However, as mentioned before, since the requirement of an HSS structure is too strong, the application of HSS matrices is limited as compared to $\mathcal{H}^2$ matrices.

**4. Error analysis.** In this section, we analyze the approximation error of SMASH. Since the HSS construction is more complicated than $\mathcal{H}^2$ construction due to the factorization of nearfield blocks, the corresponding error analysis is more involved. Here we first present error analysis for the HSS construction associated with a perfect binary tree and then an error bound for the $\mathcal{H}^2$ construction associated with a perfect $2^d$-tree can be easily derived.

For a perfect $2^d$-tree, it is natural and easy to interpret the levelwise construction in Section 3.3 in the following recursive manner:

$$A^{(l)} = U^{(l)} A^{(l-1)} \left(V^{(l)}\right)^T + B^{(l)} + E^{(l)}, \quad \forall l \geq 3, \tag{4.1}$$

where $A^{(L)} = A$ and $A^{(l)}$ ($l < L$) is a submatrix of $A$ with the following block structure:

$$(l_i, l_j) - \text{block of } A^{(l)} = \begin{cases} A|_{\bar{\mathbf{i}} \times \bar{\mathbf{j}}} & \text{if } \mathrm{lv}(i) = \mathrm{lv}(j) = l \text{ and } (i,j) \text{ is admissible}, \\ 0 & \text{if } \mathrm{lv}(i) = \mathrm{lv}(j) = l \text{ and } (i,j) \text{ is not admissible}, \end{cases} \tag{4.2}$$

$U^{(l)}, V^{(l)}, B^{(l)}$ follow the definition in Section 3.4 and $E^{(l)}$ denotes the factorization error at level $l$. The superscripts *row* and *col* used in (3.6) are dropped in (4.2) in order to simplify the notation used in the proof. Here, we assume that the items involving the index $i$ refer to rows and the items involving the index $j$ refer to columns. In addition, we introduce the notation

$$\mathbf{I}^{(l)} := \cup_{\mathrm{lv}(i)=l} \bar{\mathbf{i}} \quad \text{and} \quad \mathbf{J}^{(l)} := \cup_{\mathrm{lv}(j)=l} \bar{\mathbf{j}}.$$

Then the size of $A^{(l)}$ is equal to $|\mathbf{I}^{(l)}| \times |\mathbf{J}^{(l)}|$.

We also assume there exists a constant $r^{(l)}$ associated with each level such that

$$|\hat{\mathbf{i}}| \leq r^{(l)}, \quad |\hat{\mathbf{j}}| \leq r^{(l)}, \quad i, j \text{ at level } l, \text{ and } r^{(l+1)} \leq r^{(l)}. \tag{4.3}$$

This constant $r^{(l)}$ is actually an upper bound for the numerical ranks in the admissible blocks at level $l$.

Expanding the recursion (4.1) leads to

$$\begin{aligned}
A = & U^{(L)} \left( U^{(L-1)} \left( \cdots \left( U^{(2)} B^{(1)} \left(V^{(2)}\right)^T + B^{(2)} + E^{(2)} \right) \cdots \right) \left(V^{(L-1)}\right)^T \right. \\
& \left. + B^{(L-1)} + E^{(L-1)} \right) \left(V^{(L)}\right)^T + B^{(L)} + E^{(L)}.
\end{aligned} \tag{4.4}$$

Now we focus on the analysis of HSS approximation error. To estimate the norm of each diagonal block in $U^{(l)}$ and $V^{(l)}$, the following lemma is needed, which is a simple consequence of SRRQR in Theorem 2.4 and thus holds for both $\mathcal{H}^2$ and HSS construction.

LEMMA 4.1. *Let $s > 1$ be the prescribed elementwise bound in SRRQR in Theorem 2.4, then the norms of the $\mathcal{H}^2$ generators in (3.17) or of the HSS generators in (3.10) and (3.13) satisfy the following estimate*

$$\left\| P_i \begin{bmatrix} I \\ G_i \end{bmatrix} \right\|_F \leq s\sqrt{|\bar{\mathbf{i}}| r^{(l)}} \quad \text{and} \quad \left\| F_j \begin{bmatrix} I \\ H_j \end{bmatrix} \right\|_F \leq s\sqrt{|\bar{\mathbf{j}}| r^{(l)}}, \tag{4.5}$$

*for any nodes $i, j$ at level $l$.*

*Proof.* Under the assumption (4.3), we know that the column size of $\begin{bmatrix} I \\ G_i \end{bmatrix}$ is $|\hat{\mathbf{i}}| \leq r^{(l)}$. Since the row size of $\begin{bmatrix} I \\ G_i \end{bmatrix}$ is equal to $|\bar{\mathbf{i}}|$ and the entries of $G_i$ is bounded by $s$, we get

$$\left\| P_i \begin{bmatrix} I \\ G_i \end{bmatrix} \right\|_F \leq s\sqrt{|\hat{\mathbf{i}}| + (|\bar{\mathbf{i}}| - |\hat{\mathbf{i}}|)|\hat{\mathbf{i}}|} \leq s\sqrt{|\bar{\mathbf{i}}| r^{(l)}}.$$

The same argument applies to $\left\| F_j \begin{bmatrix} I \\ H_j \end{bmatrix} \right\|_F$. □

In the following two lemmas, Lemma 4.2 and Lemma 4.3, we investigate the local error generated in farfield approximation and nearfield approximation, respectively. Lemma 4.2 applies to both $\mathcal{H}^2$ and HSS construction, while Lemma 4.3 is only necessary for HSS construction because there is no nearfield approximation in $\mathcal{H}^2$ construction. Referring to Section 3.1, in the following we will denote by Sep the set of all well-separated pairs of subdomains corresponding to a given partitioning:

$$\text{Sep} := \{(i,j) : \Omega_i, \Omega_j \text{ are well-separated }\}. \tag{4.6}$$

LEMMA 4.2. *Suppose $i$ and $j$ are two nodes at level $l$ and $(i,j) \in \text{Sep}$. Denote the farfield approximation tolerance by $\epsilon_{far}$ as in (2.9). Then the $\mathcal{H}^2$ generators computed in (3.17) or the HSS generators computed in (3.10) and (3.13) produce the following factorization*

$$A|_{\bar{\mathbf{i}}\times\bar{\mathbf{j}}} = P_i \begin{bmatrix} I \\ G_i \end{bmatrix} A|_{\hat{\mathbf{i}}\times\hat{\mathbf{j}}} \left(F_j \begin{bmatrix} I \\ H_j \end{bmatrix}\right)^T + E_{(i,j)},$$

*where the approximation error $E_{(i,j)}$ satisfies*

$$\|E_{(i,j)}\|_F \leq 2s^2\sqrt{|\bar{\mathbf{i}}||\bar{\mathbf{j}}|}r^{(l)}\epsilon_{far}\|A|_{\bar{\mathbf{i}}\times\bar{\mathbf{j}}}\|_F.$$

*Proof.* According to (2.7) and (2.11), we know that for each $(i,j) \in \text{Sep}$

$$A|_{\bar{\mathbf{i}}\times\bar{\mathbf{j}}} = \hat{U}_{\bar{\mathbf{i}}}\hat{B}_{i,j}\hat{V}_{\bar{\mathbf{j}}}^T + E_F|_{\bar{\mathbf{i}}\times\bar{\mathbf{j}}} = P_i \begin{bmatrix} I \\ G_i \end{bmatrix} \hat{U}_{\hat{\mathbf{i}}}\hat{B}_{i,j}\hat{V}_{\hat{\mathbf{j}}}^T \left(F_j \begin{bmatrix} I \\ H_j \end{bmatrix}\right)^T + E_F|_{\bar{\mathbf{i}}\times\bar{\mathbf{j}}}$$

$$= P_i \begin{bmatrix} I \\ G_i \end{bmatrix} A|_{\hat{\mathbf{i}}\times\hat{\mathbf{j}}} \left(F_j \begin{bmatrix} I \\ H_j \end{bmatrix}\right)^T + E_{(i,j)},$$

where

$$E_{(i,j)} := -P_i \begin{bmatrix} I \\ G_i \end{bmatrix} E_F|_{\hat{\mathbf{i}}\times\hat{\mathbf{j}}} \left(F_j \begin{bmatrix} I \\ H_j \end{bmatrix}\right)^T + E_F|_{\bar{\mathbf{i}}\times\bar{\mathbf{j}}}. \tag{4.7}$$

Based on (2.9) and Lemma 4.1, we deduce that

$$\|E_{(i,j)}\|_F \leq \|P_i \begin{bmatrix} I \\ G_i \end{bmatrix}\|_F \|E_F|_{\hat{\mathbf{i}}\times\hat{\mathbf{j}}}\|_F \|F_j \begin{bmatrix} I \\ H_j \end{bmatrix}\|_F + \|E_F|_{\bar{\mathbf{i}}\times\bar{\mathbf{j}}}\|_F$$

$$\leq 2s\sqrt{|\bar{\mathbf{i}}|r^{(l)}}\epsilon_{\text{far}}\|A|_{\bar{\mathbf{i}}\times\bar{\mathbf{j}}}\|_F s\sqrt{|\bar{\mathbf{j}}|r^{(l)}} = 2s^2\sqrt{|\bar{\mathbf{i}}||\bar{\mathbf{j}}|}r^{(l)}\epsilon_{\text{far}}\|A|_{\bar{\mathbf{i}}\times\bar{\mathbf{j}}}\|_F.$$

☐

LEMMA 4.3. *Suppose $i$ and $j$ are two nodes at level $l$ of the HSS tree and $(i,j) \notin \text{Sep}$. Then in HSS construction, the application of the truncated SVD with relative tolerance $\epsilon_{SVD}$ in (3.7) and (3.11) produces the following factorization*

$$A|_{\bar{\mathbf{i}}\times\bar{\mathbf{j}}} = S_i \Sigma_i^- \widetilde{S}_j + E_{(i,j)},$$

*where*

$$\|E_{(i,j)}\|_F \leq s^2|\bar{\mathbf{j}}|\sqrt{|\bar{\mathbf{i}}|r^{(l)}}r^{(l)}\epsilon_{SVD}\|A_j^|\|_F + s|\bar{\mathbf{i}}|\sqrt{|\bar{\mathbf{j}}|r^{(l)}}\epsilon_{SVD}\|A_i^-\|_F.$$

*Proof.* From (3.7) and (3.11), we know that

$$A|_{\bar{\mathbf{i}}\times\bar{\mathbf{j}}} = S_i \Sigma_i^- \widetilde{S}_j + E_\Sigma^-|_{\bar{\mathbf{i}}\times\bar{\mathbf{j}}} = P_i \begin{bmatrix} I \\ G_i \end{bmatrix} (A|_{\hat{\mathbf{i}}\times\bar{\mathbf{j}}} - E_\Sigma^-|_{\hat{\mathbf{i}}\times\bar{\mathbf{j}}}) + E_\Sigma^-|_{\bar{\mathbf{i}}\times\bar{\mathbf{j}}}$$

$$= P_i \begin{bmatrix} I \\ G_i \end{bmatrix} A|_{\hat{\mathbf{i}}\times\hat{\mathbf{j}}} \left( F_j \begin{bmatrix} I \\ H_j \end{bmatrix} \right)^T + E_{(i,j)},$$

where

(4.8) $$E_{(i,j)} := - P_i \begin{bmatrix} I \\ G_i \end{bmatrix} E_\Sigma^|_{\hat{\mathbf{i}}\times\hat{\mathbf{j}}} \left( F_j \begin{bmatrix} I \\ H_j \end{bmatrix} \right)^T + P_i \begin{bmatrix} I \\ G_i \end{bmatrix} \left( E_\Sigma^|_{\hat{\mathbf{i}}\times\bar{\mathbf{j}}} - E_\Sigma^-|_{\hat{\mathbf{i}}\times\bar{\mathbf{j}}} \right) + E_\Sigma^-|_{\bar{\mathbf{i}}\times\bar{\mathbf{j}}}.$$

Introduce the notation

$$\check{\mathbf{i}} := \bar{\mathbf{i}} \setminus \hat{\mathbf{i}} \quad \text{and} \quad \check{\mathbf{j}} := \bar{\mathbf{j}} \setminus \hat{\mathbf{j}},$$

we then have

$$E_\Sigma^|_{\hat{\mathbf{i}}\times\bar{\mathbf{j}}} = E_\Sigma^|_{\hat{\mathbf{i}}\times\hat{\mathbf{j}}} \left( F_j \begin{bmatrix} I \\ 0 \end{bmatrix} \right) + E_\Sigma^|_{\hat{\mathbf{i}}\times\check{\mathbf{j}}} \left( F_j \begin{bmatrix} 0 \\ I \end{bmatrix} \right) \quad \text{and} \quad E_\Sigma^-|_{\bar{\mathbf{i}}\times\bar{\mathbf{j}}} = P_i \begin{bmatrix} I \\ 0 \end{bmatrix} E_\Sigma^-|_{\hat{\mathbf{i}}\times\bar{\mathbf{j}}} + P_i \begin{bmatrix} 0 \\ I \end{bmatrix} E_\Sigma^-|_{\check{\mathbf{i}}\times\bar{\mathbf{j}}}.$$

Substituting the above identities into (4.8), we obtain

(4.9) $$E_{(i,j)} = - P_i \begin{bmatrix} I \\ G_i \end{bmatrix} E_\Sigma^|_{\hat{\mathbf{i}}\times\hat{\mathbf{j}}} \left( F_j \begin{bmatrix} 0 \\ H_j \end{bmatrix} \right)^T + P_i \begin{bmatrix} I \\ G_i \end{bmatrix} E_\Sigma^|_{\hat{\mathbf{i}}\times\check{\mathbf{j}}} \left( F_j \begin{bmatrix} 0 \\ I \end{bmatrix} \right)^T$$
$$- P_i \begin{bmatrix} 0 \\ G_i \end{bmatrix} E_\Sigma^-|_{\hat{\mathbf{i}}\times\bar{\mathbf{j}}} + P_i \begin{bmatrix} 0 \\ I \end{bmatrix} E_\Sigma^-|_{\check{\mathbf{i}}\times\bar{\mathbf{j}}}.$$

Based on (3.8), (3.12) and Lemma 4.1, we have the estimate

$$\|E_{(i,j)}\|_F \leq \left\| \begin{bmatrix} I \\ G_i \end{bmatrix} \right\|_F \left( \|H_j\|_F \|E_\Sigma^|_{\hat{\mathbf{i}}\times\hat{\mathbf{j}}}\|_F + \|E_\Sigma^|_{\hat{\mathbf{i}}\times\check{\mathbf{j}}}\|_F \right) + \|G_i\|_F \|E_\Sigma^-|_{\hat{\mathbf{i}}\times\bar{\mathbf{j}}}\|_F + \|E_\Sigma^-|_{\check{\mathbf{i}}\times\bar{\mathbf{j}}}\|_F$$

$$\leq \left\| \begin{bmatrix} I \\ G_i \end{bmatrix} \right\|_F \sqrt{\|H_j\|_F^2 + 1} \|E_\Sigma^|_{\hat{\mathbf{i}}\times\bar{\mathbf{j}}}\|_F + \sqrt{\|G_i\|_F^2 + 1} \|E_\Sigma^-|_{\bar{\mathbf{i}}\times\bar{\mathbf{j}}}\|_F$$

$$\leq s^2 |\bar{\mathbf{j}}| \sqrt{|\bar{\mathbf{i}}|r^{(l)}} r^{(l)} \epsilon_{\text{SVD}} \|A_j^|\|_F + s|\bar{\mathbf{i}}| \sqrt{|\bar{\mathbf{j}}|r^{(l)}} \epsilon_{\text{SVD}} \|A_i^-\|_F.$$

☐

Based on Lemmas 4.2–4.3, we can estimate the total approximation error at level $l$ in HSS construction.

LEMMA 4.4. *Assume $\epsilon_{SVD}$ and $\epsilon_{far}$ are the approximation tolerances used in the approximation of nearfield and farfield blocks in HSS construction. Then the approximation error $E^{(l)}$ in (4.1) satisfies the following bound*

$$\|E^{(l)}\|_F \leq s^2 2^{l/2+2} \left( r^{(l+1)} \right)^{3/2} \left( r^{(l)} \right)^{3/2} \epsilon_{SVD} \|A^{(l)}\|_F + 4s^2 r^{(l)} r^{(l+1)} \epsilon_{far} \|A^{(l)}\|_F,$$

*where we have set $r^{(L+1)} := r^{(L)}$.*

*Proof.* Based on Lemma 4.3, we know that the approximation error from the nearfield compression at level $l$ can be estimated as follows:

$$\sum_{(i,j)\notin \text{Sep}} \|E_{(i,j)}\|_F^2 \leq 2 \sum_{(i,j)\notin \text{Sep}} \left( s^4 |\bar{\mathbf{j}}|^2 |\bar{\mathbf{i}}| \left( r^{(l)} \right)^3 \epsilon_{\text{SVD}}^2 \|A_j^|\|_F^2 + s^2 |\bar{\mathbf{i}}|^2 |\bar{\mathbf{j}}| r^{(l)} \epsilon_{\text{SVD}}^2 \|A_i^-\|_F^2 \right)$$

$$\leq 2 \sum_{i,j \text{ at level } l} \left( s^4 |\bar{\mathbf{j}}|^2 |\bar{\mathbf{i}}| \left( r^{(l)} \right)^3 \epsilon_{\text{SVD}}^2 \|A_j^|\|_F^2 + s^2 |\bar{\mathbf{i}}|^2 |\bar{\mathbf{j}}| r^{(l)} \epsilon_{\text{SVD}}^2 \|A_i^-\|_F^2 \right).$$

Since $|\bar{\mathbf{i}}| \leq 2r^{(l+1)}, |\bar{\mathbf{j}}| \leq 2r^{(l+1)}$ for nodes $i, j$ at level $l$ and there are $2^{l-1}$ nodes at this level, we further have

$$\sum_{(i,j)\notin \text{Sep}} \|E_{(i,j)}\|_F^2 \leq \sum_{i,j \text{ at level } l} \left(2s^4 \left(2r^{(l+1)}\right)^3 \left(r^{(l)}\right)^3 \epsilon_{\text{SVD}}^2 \|A_j^{|}\|_F^2 + 2s^2 \left(2r^{(l+1)}\right)^3 r^{(l)} \epsilon_{\text{SVD}}^2 \|A_i^{-}\|_F^2\right)$$

$$\leq \sum_{lv(i)=l} 2s^4 \left(2r^{(l+1)}\right)^3 \left(r^{(l)}\right)^3 \epsilon_{\text{SVD}}^2 \|A^{(l)}\|_F^2$$

$$+ \sum_{lv(j)=l} 2s^2 \left(2r^{(l+1)}\right)^3 r^{(l)} \epsilon_{\text{SVD}}^2 \|A^{(l)}\|_F^2$$

$$\leq 2^l \left(s^4 \left(2r^{(l+1)}\right)^3 \left(r^{(l)}\right)^3 + s^2 \left(2r^{(l+1)}\right)^3 r^{(l)}\right) \epsilon_{\text{SVD}}^2 \|A^{(l)}\|_F^2$$

$$\leq s^4 2^{l+4} \left(r^{(l+1)}\right)^3 \left(r^{(l)}\right)^3 \epsilon_{\text{SVD}}^2 \|A^{(l)}\|_F^2.$$

The approximation error for the farfield compression at level $l$ can be estimated as follows:

$$\sum_{(i,j)\in \text{Sep}} \|E_{(i,j)}\|_F^2 \leq \sum_{i,j \text{ at level } l} 4s^4 \left(2r^{(l+1)} r^{(l)}\right)^2 \epsilon_{\text{far}}^2 \|A|_{\bar{\mathbf{i}} \times \bar{\mathbf{j}}}\|_F^2 = 16s^4 \left(r^{(l+1)} r^{(l)}\right)^2 \epsilon_{\text{far}}^2 \|A^{(l)}\|_F^2.$$

To sum up both farfield and nearfield approximation errors, we obtain the estimate for the overall approximation error introduced at level $l$:

$$\|E^{(l)}\|_F = \left(\sum_{(i,j)\notin \text{Sep}} \|E_{(i,j)}\|_F^2 + \sum_{(i,j)\in \text{Sep}} \|E_{(i,j)}\|_F^2\right)^{\frac{1}{2}}$$

$$\leq s^2 2^{l/2+2} \left(r^{(l+1)}\right)^{3/2} \left(r^{(l)}\right)^{3/2} \epsilon_{\text{SVD}} \|A^{(l)}\|_F + 4s^2 r^{(l)} r^{(l+1)} \epsilon_{\text{far}} \|A^{(l)}\|_F.$$

☐

The overall HSS approximation error in the Frobenius norm can then be derived in the following theorem.

THEOREM 4.5. *Suppose the HSS tree has L levels. With the assumptions in Lemma 4.4, SMASH produces the following factorization*

$$(4.10) \quad \begin{aligned} A = & U^{(L)} \left(U^{(L-1)} \left(\ldots \left(U^{(2)} B^{(1)} \left(V^{(2)}\right)^T + B^{(2)}\right) \ldots \right) \left(V^{(L-1)}\right)^T \\ & + B^{(L-1)}\right) \left(V^{(L)}\right)^T + B^{(L)} + E, \end{aligned}$$

*where the approximation error $E$ satisfies the estimate*

$$\|E\|_F \leq C_1 \epsilon_{SVD} \|A\|_F + C_2 \epsilon_{far} \|A\|_F,$$

*with*

$$C_1 = \sum_{l=2}^{L-1} 2^{L+l/2+2} s^{2L-2l+2} \left(r^{(l+1)} \ldots r^{(L)}\right)^2 \left(r^{(l+1)}\right)^{3/2} \left(r^{(l)}\right)^{5/2}$$

$$C_2 = \sum_{l=2}^{L} 2^{L+2} s^{2L-2l+2} \left(r^{(l)} r^{(l+1)} \ldots r^{(L)}\right)^2 r^{(l+1)}.$$

*Proof.* According to (4.4), we know that the overall HSS approximation error $E$ has the expression

$$(4.11) \quad E = \left(U^{(L)}\ldots U^{(3)}\right) E^{(2)} \left(V^{(L)}\ldots V^{(3)}\right)^T + \cdots + U^{(L)} E^{(L-1)} \left(V^{(L)}\right)^T + E^{(L)}.$$

Note that the column size of $\left(U^{(L)}\ldots U^{(l+1)}\right)$ is bounded by $r^{(l)}2^l$. Thus we have

$$(4.12) \quad \|U^{(L)}\ldots U^{(l+1)}\|_F \leq \sqrt{r^{(l)}2^l}\|U^{(L)}\ldots U^{(l+1)}\|_2 \leq 2^{L/2} s^{L-l}\sqrt{r^{(l)}} r^{(l+1)}\ldots r^{(L)},$$

and same upper bound holds for $\|V^{(L)}\ldots V^{(l+1)}\|_F$. It follows from (4.11) that

$$\|E\|_F \leq \sum_{l=2}^{L-1} \|U^{(L)}\ldots U^{(l+1)}\|_F \|E^{(l)}\|_F \|V^{(L)}\ldots V^{(l+1)}\|_F + \|E^{(L)}\|_F$$

$$\leq \sum_{l=2}^{L-1} 2^L s^{2L-2l} r^{(l)} \left(r^{(l+1)}\ldots r^{(L)}\right)^2 s^2 2^{l/2+2} \left(r^{(l+1)}\right)^{3/2} \left(r^{(l)}\right)^{3/2} \epsilon_{\text{SVD}} \|A^{(l)}\|_F$$

$$(4.13) \quad + \sum_{l=2}^{L-1} 2^L s^{2L-2l} r^{(l)} \left(r^{(l+1)}\ldots r^{(L)}\right)^2 4 s^2 r^{(l+1)} r^{(l)} \epsilon_{\text{far}} \|A^{(l)}\|_F + \|E^{(L)}\|_F$$

$$\leq \sum_{l=2}^{L-1} 2^{L+l/2+2} s^{2L-2l+2} \left(r^{(l+1)}\ldots r^{(L)}\right)^2 \left(r^{(l+1)}\right)^{3/2} \left(r^{(l)}\right)^{5/2} \epsilon_{\text{SVD}} \|A\|_F$$

$$+ \sum_{l=2}^{L} 2^{L+2} s^{2L-2l+2} \left(r^{(l)} r^{(l+1)}\ldots r^{(L)}\right)^2 r^{(l+1)} \epsilon_{\text{far}} \|A\|_F.$$

□

It is easy to see that the bound in (4.13) is quite pessimistic because we bound $\|A^{(l)}\|_F$ from above by $\|A\|_F$, where $A^{(l)}$ (defined in (4.2)) is a submatrix of $A$ of size $|\mathbf{I}^{(l)}| \times |\mathbf{J}^{(l)}|$.

COROLLARY 4.6. *Besides the assumptions in Theorem 4.5, if there also exists a constant $r \geq 2$ such that $r^{(l)} \leq r$ for each $l > 1$, then the approximation error $E$ in SMASH satisfies the estimate:*

$$\|E\|_F \leq (2r^2 s^2)^L (16\epsilon_{SVD} + 8\epsilon_{far}) \|A\|_F.$$

Next we estimate the error in $\mathcal{H}^2$ approximation. Note that, by setting $\epsilon_{\text{SVD}} = 0$, the error estimate in Theorem 4.5 also holds for the $\mathcal{H}^2$ construction with a possibly different constant $C_2$ depending on dimension $d$. In fact, analogous to Lemma 4.4, we first have the following estimate for the $\mathcal{H}^2$ construction:

$$\|E^{(l)}\|_F \leq 4s^2 r^{(l)} r^{(l+1)} \epsilon_{\text{far}} \|A^{(l)}\|_F,$$

where $E^{(l)}$ denotes the approximation error introduced at level $l$, as defined in (4.1). Assume the $\mathcal{H}^2$ matrix is associated with a perfect $2^d$-tree $\mathcal{T}$ ($d \in \{1,2,3\}$), i.e., each nonleaf node of $\mathcal{T}$ has $2^d$ children and all leaves are at the same level. For a node $i$ at level $l$, since $|\bar{\mathbf{i}}| \leq 2^d r^{(l+1)} \leq 2^d r^{(l)}$, it follows from (3.22), (3.17) and Lemma 4.1 that

$$\|U^{(l)}\|_2 \leq s 2^{d/2} r^{(l)}.$$

Notice that the column size of $\left(U^{(L)}\ldots U^{(l+1)}\right)$ is bounded by $r^{(l)}2^{dl}$. Then the counterpart of (4.12) can be obtained for the $\mathcal{H}^2$ construction as:

$$\|U^{(L)}\ldots U^{(l+1)}\|_F \leq \sqrt{r^{(l)}2^{dl}} \|U^{(L)}\ldots U^{(l+1)}\|_2 \leq 2^{dL/2} s^{L-l} \sqrt{r^{(l)}} r^{(l+1)}\ldots r^{(L)},$$

and same upper bound holds for $\|V^{(L)} \ldots V^{(l+1)}\|_F$. The total approximation error can now be obtained by means of (4.11):

$$\begin{aligned}
\|E\|_F &\leq \sum_{l=2}^{L-1} \|U^{(L)} \ldots U^{(l+1)}\|_F \|E^{(l)}\|_F \|V^{(L)} \ldots V^{(l+1)}\|_F + \|E^{(L)}\|_F \\
&\leq \sum_{l=2}^{L-1} 2^{dL} s^{2L-2l} r^{(l)} \left(r^{(l+1)} \ldots r^{(L)}\right)^2 4s^2 r^{(l+1)} r^{(l)} \epsilon_{\text{far}} \|A^{(l)}\|_F + 4s^2 (r^{(L)})^2 \epsilon_{\text{far}} \|A\|_F \\
&\leq \sum_{l=2}^{L-1} 2^{dL+2} s^{2L-2l+2} \left(r^{(l)} r^{(l+1)} \ldots r^{(L)}\right)^2 r^{(l+1)} \epsilon_{\text{far}} \|A\|_F + 4s^2 (r^{(L)})^2 \epsilon_{\text{far}} \|A\|_F \\
&\leq \sum_{l=2}^{L} 2^{dL+2} s^{2L-2l+2} \left(r^{(l)} r^{(l+1)} \ldots r^{(L)}\right)^2 r^{(l+1)} \epsilon_{\text{far}} \|A\|_F.
\end{aligned}$$

The above error analysis yields the following theorem.

THEOREM 4.7. *Suppose $\mathcal{T}$ is a perfect $2^d$-tree with $L$ levels, associated with the $\mathcal{H}^2$ approximation of $A$. Under the assumptions in (4.3), SMASH produces the following factorization*

$$\begin{aligned}
(4.14) \quad A = &U^{(L)} \left( U^{(L-1)} \left( \ldots \left( U^{(2)} B^{(1)} \left( V^{(2)} \right)^T + B^{(2)} \right) \ldots \right) \left( V^{(L-1)} \right)^T \right. \\
&\left. + B^{(L-1)} \right) \left( V^{(L)} \right)^T + B^{(L)} + E,
\end{aligned}$$

*where the approximation error $E$ satisfies the estimate*

$$\|E\|_F \leq C \epsilon_{far} \|A\|_F,$$

*with*

$$C = \sum_{l=2}^{L} 2^{dL+2} s^{2L-2l+2} \left(r^{(l)} r^{(l+1)} \ldots r^{(L)}\right)^2 r^{(l+1)}.$$

COROLLARY 4.8. *Besides the assumptions in Theorem 4.7, if there also exists a constant $r \geq 2$ such that $r^{(l)} \leq r$ for each $l > 1$, then the approximation error $E$ in SMASH satisfies the estimate:*

$$\|E\|_F \leq (2^d r^2 s^2)^L 8 \epsilon_{far} \|A\|_F.$$

**5. Complexity analysis.** This section studies the complexity of SMASH for an $n \times n$ matrix. For simplicity, we only consider the case when $X = Y$ and the points are uniformly distributed. Under this assumption, a perfect tree $\mathcal{T}$ will be used for both HSS and $\mathcal{H}^2$ structures.

**5.1. Complexity for the HSS construction.** We start with the HSS construction case. Suppose $\mathcal{T}$ has $L$ levels such that $n = O(r 2^L)$, where $r$ is a positive integer such that the rank of HSS generators is bounded above by $r$ and

(5.1) $$|\bar{\mathbf{i}}^{row}| \leq 2r, \quad |\bar{\mathbf{i}}^{col}| \leq 2r, \quad \forall i \neq \text{root}.$$

Notice that in the context of integral equations in potential theory, the assumption (5.1) in general holds only for integral equations defined on a curve. Since the points are uniformly distributed, for each nonroot node $i$, the number of nodes in $\mathcal{N}_i$ is very small, which we assume to be bounded above by 3. Under these assumptions, we have the following complexity estimate.

THEOREM 5.1. *Let $\mathcal{T}$ be a perfect binary tree with $L$ levels and* (5.1) *hold. Then the complexity of SMASH for the HSS construction in Section 3.3.1 is $O(n)$.*

*Proof.* Based on (5.1), it is easy to see that, for each nonroot node $i$, the compression cost for its nearfield blocks in (3.7) is $O(r^3)$. This is because the size of the nearfield block row in (3.7) is no larger than $2r$-by-$6r$ under the above assumption for $\mathcal{N}_i$. Besides, the farfield basis matrix $\hat{U}_i$ has column size at most $r$, so the cost of an SRRQR procedure in (3.9) is $O(r^3)$. Therefore, the compression cost associated with each nonroot node $i$ is $O(r^3)$ and the complexity of the HSS construction is $O(2^L r^3) = r^2 O(n) = O(n)$. □

**5.2. Complexity for the $\mathcal{H}^2$ construction.** For the $\mathcal{H}^2$ construction case, we assume that when $X \subset \mathbb{R}^d$, $\mathcal{T}$ is a perfect $2^d$-tree with $L$ levels such that $n = O(r2^{dL})$ and $r$ is a positive integer such that the rank of $\mathcal{H}^2$ generators is bounded by $r$ and

$$(5.2) \qquad |\bar{\mathbf{i}}^{row}| \leq r2^d, \quad |\bar{\mathbf{i}}^{col}| \leq r2^d, \quad \forall i \neq \text{root}.$$

The analysis here is simpler than that of the HSS construction in Section 5.1. Since each node $i$ only involves the compression of farfield basis $\hat{U}_i|_{\bar{\mathbf{i}}^{row}}$ (as well as $\hat{V}_i|_{\bar{\mathbf{i}}^{col}}$), whose size is no larger than $r2^d$-by-$r$ under the assumption (5.2), we deduce that the compression cost associated with each node is $O(r^3)$. As a result, the complexity of the $\mathcal{H}^2$ construction is $O(2^{dL} r^3) = r^2 O(n)$. Thus we conclude:

THEOREM 5.2. *Let $\mathcal{T}$ be a perfect $2^d$-tree with $L$ levels and* (5.2) *hold. Then the complexity of SMASH for the $\mathcal{H}^2$ construction in Section 3.3.2 is $O(n)$.*

**6. Numerical examples.** In this section, we present numerical examples to illustrate the performance of SMASH. All of the numerical results were performed in MATLAB R2014b on a macbook air with a 1.6 GHz CPU and 8 GB of RAM. The following notation is used throughout the section:
- $n$: the size of $A$;
- $t_{\text{constr}}$: wall clock time for constructing $\hat{A}$ in seconds;
- $t_{\text{matvec}}$: wall clock time for multiplying $\hat{A}$ with a vector in seconds;
- $t_{\text{sol}}$: wall clock time for solving $\hat{A}x = b$ in seconds;
- $\epsilon_{svd}$: relative tolerance used in the truncated SVD for the nearfield compression;
- rand$([0,1])$: a random number sampled from the uniform distribution in $[0,1]$.

**6.1. Choice of parameters.** Since the quality of a degenerate approximation depends on the underlying kernel function, there is no rule of thumb in general on choosing the parameters to satisfy a prescribed tolerance. For completeness, here we present a heuristic approach that we use in all numerical experiments on the choice of parameters.

Given a matrix $A$ and a tolerance $\epsilon$, suppose one wants to construct a hierarchical matrix $\hat{A}$ ($\mathcal{H}$, $\mathcal{H}^2$, or HSS) such that $\|A - \hat{A}\|_{\max} \approx \epsilon$. Then the following approach is adopted to determine parameters $\tau, r$.

The choice of separation ratio $\tau \in (0,1)$ only depends on the dimension of the problem, so it is chosen first. We choose $\tau$ such that $\tau \leq 0.7$ and, in general, a slightly larger $\tau$ is preferred for higher dimensional problems. For example, we choose $\tau = 0.6$ for essentially one-dimensional problems, such as those in Section 6.3 and Section 6.4; we choose $\tau = 0.65$ for two-dimensional problems, as in Section 6.2.

Having chosen a separation ratio $\tau$, we use the following function to determine the farfield

Table 6.1: Construction and matrix-vector multiplication of an $\mathcal{H}^2$ matrix.

| $n = m^2$ | $\|\hat{A}u - Au\|/\|Au\|$ | $t_{\text{constr}}$ | $t_{\text{matvec}}$ |
|---|---|---|---|
| 1600 | $6.69 \times 10^{-13}$ | 0.52 | 0.02 |
| 6400 | $2.00 \times 10^{-12}$ | 1.97 | 0.07 |
| 25600 | $3.65 \times 10^{-12}$ | 9.53 | 0.30 |
| 102400 | $4.87 \times 10^{-12}$ | 39.47 | 1.18 |

approximation rank $r$ used in constructing $\hat{U}_i, \hat{V}_i$ (before the SRRQR postprocessing):

$$r = \begin{cases} \lfloor \log \epsilon / \log \tau - 20 \rfloor, & \text{if } \epsilon < 10^{-8}, \\ \lfloor \log \epsilon / \log \tau - 15 \rfloor, & \text{if } 10^{-8} \leq \epsilon < 10^{-6}, \\ \max\{\lfloor \log \epsilon / \log \tau - 10 \rfloor, 5\} & \text{otherwise}, \end{cases}$$

where $\lfloor x \rfloor$ yields the largest integer less than or equal to $x$. For example, in Section 6.2, $\epsilon = 10^{-7}$, $\tau = 0.65$, $r = 22$; in Section 6.3, $\epsilon = 10^{-8}$, $\tau = 0.6, r = 21$; in Section 6.4, $\epsilon = 10^{-10}$, $\tau = 0.6$, $r = 25$.

**6.2. Construction and matrix-vector multiplication of $\mathcal{H}^2$ matrices.** We consider the construction and matrix-vector multiplication of an $\mathcal{H}^2$ approximation associated with the kernel in (2.1) with $d_x = 1$. We chose $X$ as a uniform $m \times m$ grid in $[0, 1]^2$ and set $Y = X$. The computational domain $[0, 1]^2$ was recursively divided into 4 subdomains until the number of the points inside each domain was less than or equal to 50. We used the truncated Taylor expansion (2.4) with $r = 22$ terms and the separation ratio $\tau = 0.65$ to compress farfield blocks. The approximation error of this construction is measured by the relative error $\|\hat{A}u - Au\|/\|Au\|$, where $u$ is a random vector of length $n = m^2$ with entries generated by rand$([0, 1])$. The numerical results are reported in Table 6.1.

As can be seen from Table 6.1, SMASH for the $\mathcal{H}^2$ construction and the matrix-vector multiplication described in Section 3.4 scale linearly, which is consistent with the complexity analysis in Section 5.

**6.3. Cauchy-like matrices.** We consider in this section the numerical solution of Cauchy-like matrices. It is known that Cauchy-like matrices are related to other types of structured matrices including Toeplitz matrices, Vandermonde matrices, Hankel matrices and their variants [45, 47, 46, 48]. Consider the kernel $\kappa(x, y) = 1/(x - y)$, $x \neq y \in \mathbb{C}$. Let $x_i, y_j (i, j = 1 : n)$ be $2n$ pairwise distinct points in $\mathbb{C}$. The Cauchy matrix is then given by $\mathcal{C} = [\kappa(x_i, y_j)]_{i,j=1:n}$, which is known to be invertible [19]. Given two matrices $w, v \in \mathbb{C}^{n \times p}$, the $(i, j)$-entry of a Cauchy-like matrix $A$ associated with generators $w, v$ is defined by [9]

$$(6.1) \qquad a_{i,j} = \frac{1}{x_i - y_j} \sum_{l=1}^{p} w_{i,l} v_{j,l}.$$

For simplicity, we consider the case $p = 2$, i.e., $w$ (as well as $v$) is composed of two column vectors. Denote by $\hat{w}_1, \hat{w}_2, \hat{v}_1, \hat{v}_2$ the column vectors in $u, v$, i.e., $w = [\hat{w}_1, \hat{w}_2], v = [\hat{v}_1, \hat{v}_2]$. It can be seen that $A$ can be written as

$$(6.2) \qquad A = \text{diag}(\hat{w}_1)\mathcal{C}\text{diag}(\hat{v}_1) + \text{diag}(\hat{w}_2)\mathcal{C}\text{diag}(\hat{v}_2).$$

Existing approaches for solving Cauchy or Cauchy-like linear systems associated with points in $\mathbb{R}$ mainly rely on some variants of Gaussian elimination with pivoting techniques. For example, fast $O(n^2)$ algorithms for solving Cauchy linear systems can be found in [16, 17, 24, 29], etc.; a superfast $O(n \log^3 n)$ algorithm based on a sequential block Gaussian elimination process was proposed in [49]. The performance of most existing methods depends on the the distribution of point sets $x, y$. As pointed out in [17], if two sets of points $x, y$ can not be separated, for example, when they are interlaced, existing algorithms (for example, BP-type algorithm of [16]) suffer from backward stability issues. Moreover, due to the use of pivoting techniques, the accuracy of existing algorithms heavily depend on the ordering of points [16, 17] and the analysis is limited to the case when the points are in $\mathbb{R}$.

Therefore, in view of the issues mentioned above, we assume $x_i, y_j$ are mixed together such that in adaptive partitioning (see Section 3.1), each box contains the same number of points from $x_i$ and $y_j$. We also consider that $x_i, y_j$ are distributed on a curve in $\mathbb{R}^2$ as illustrated in Fig. 6.2 to demonstrate that the algorithm is independent of the ordering of points and is applicable for points in $\mathbb{C}$.

We construct the HSS approximation $\hat{A}$ to $A$ using SMASH discussed in Section 3.3.1 and then solve the linear system associated with $\hat{A}$ using a fast $ULV$ factorization solver [22]. Due to the choice of stable expansion in (2.4), arbitrarily high approximation accuracy can be achieved without stability issues [18].

Note that the HSS approximation to $\mathcal{C}$ can be readily obtained as in Section 3. Consequently, the HSS approximations to $\text{diag}(\hat{w}_1)\mathcal{C}\text{diag}(\hat{v}_1)$ and $\text{diag}(\hat{w}_2)\mathcal{C}\text{diag}(\hat{v}_2)$ can be derived by modifying $U, V, D$ generators, respectively. The sum of these two HSS representations is also an HSS matrix whose generators can be easily obtained using the technique presented in [57] by merging the two sets of HSS generators. Hence the HSS approximation to $A$ is derived.

In the first experiment, the point sets $\{x_k\}_{k=1}^n, \{y_k\}_{k=1}^n$ are chosen as follows:

(6.3) $$x_k = k/(n+1), \quad y_k = x_k + 10^{-7} * \text{rand}([0,1]), \quad k = 1, \ldots, n.$$

In the second experiment, the point sets are distributed on the curve illustrated in Fig. 6.1 that is parametrized by

$$\gamma(t) = e^{-\pi i/6} * [(0.5 + \sin(4\pi t))\cos(2\pi t) + i(0.5 + \sin(4\pi t))\sin(2\pi t)], \quad t \in [0,1],$$

and $\{x_k\}_{k=1}^n, \{y_k\}_{k=1}^n$ are given by

(6.4) $$x_k = \gamma(k/(n+1)), \quad y_k = \gamma(x_k + 10^{-7} * \text{rand}([0,1])), \quad k = 1, \ldots, n.$$

In the third experiment, the point set $\{x_k\}_{k=1}^n$ is on the snail geometry in $\mathbb{C}$ as illustrated in Fig. 6.2 and the point set $\{y_k\}_{k=1}^n$ is given by

(6.5) $$y_k = x_k + 10^{-7} * \text{rand}([0,1]), \quad k = 1, \ldots, n.$$

For generators of this Cauchy-like matrix, i.e., $w = [\hat{w}_1, \hat{w}_2], v = [\hat{v}_1, \hat{v}_2]$, we chose $\hat{w}_l, \hat{v}_l (l = 1, 2)$ such that each entry in those vectors was given by $\text{rand}([0,1])$. In order to solve the linear system $Au = b$, we constructed an HSS approximation $\hat{A}$ to $A$ in (6.2) by SMASH. 1D boxes (i.e., intervals) and 2D boxes (i.e., rectangles) are used in adaptive partitioning for point sets in (6.3) and (6.5), respectively. The right subfigures in Fig. 6.2 illustrate the adaptive partitioning using 2D boxes as described in Remark A.1. We chose separation ratio as $\tau = 0.6$ and adaptive partitioning stopped when the number of points inside each box was less than or equal 50. The nearfield blocks were compressed through SVD with truncation tolerance $10^{-9}$. The exact solution was set to be a column vector $u$ of length $n$ with entries generated by $\text{rand}([0,1])$, and the right-hand side $b$ was formed by $b = Au$.

The numerical results for three Cauchy-like matrix problems are reported in Table 6.2. From Table 6.2, we see that the computational time for both construction and the solution scale linearly, and SMASH in Section 3.3.1 is quite robust with respect to complex geometries. Moreover, it can be seen that SMASH is independent of the ordering of points.

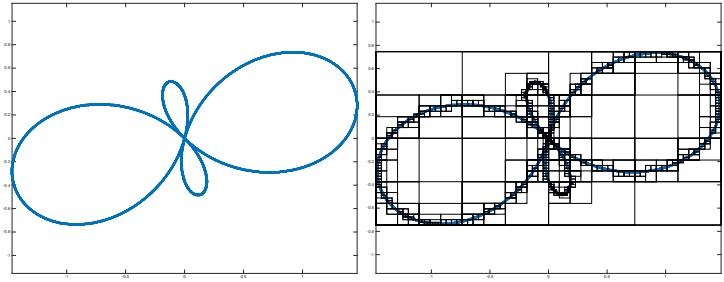

Fig. 6.1: Honeybee geometry used for the numerical experiments in Table 6.2. Left: Original curve; Right: Adaptive partitioning of the curve for the case when $n = 12800$ in Table 6.2.

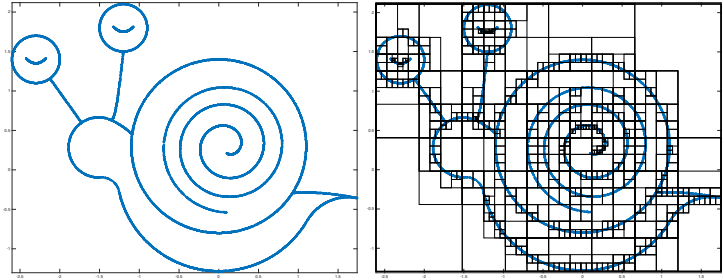

Fig. 6.2: Snail geometry used for the numerical experiments in Table 6.2. Left: Original curve; Right: Adaptive partitioning of the curve for the case when $n = 12800$ in Table 6.2.

**6.4. Integral equations.** In this section, we solve Laplace boundary value problems via the integral equation method. Assume $\Omega$ is a smooth simply-connected domain in $\mathbb{R}^2$ and let $\Gamma = \partial \Omega$ be the boundary of $\Omega$ of class $C^2$. Consider the interior Dirichlet problem: find $u \in C^2(\Omega) \cap C(\overline{\Omega})$ such that

$$\begin{aligned} \Delta u &= 0 \quad \text{in } \Omega, \\ u &= u_D \quad \text{on } \Gamma, \end{aligned} \tag{6.6}$$

where $u_D \in C(\Gamma)$ is given. With smooth boundary curves and Dirichlet data, the wellposedness of this problem is well studied in potential theory [38, 39].

The fundamental solution and its gradient (in terms of $\mathbf{y}$) for the Laplace equation in $\mathbb{R}^2$ are given by:

$$\Phi(\mathbf{x}, \mathbf{y}) = -\frac{1}{2\pi} \log |\mathbf{x} - \mathbf{y}|, \quad \text{and} \quad \nabla_y \Phi(\mathbf{x}, \mathbf{y}) = -\frac{1}{2\pi} \frac{\mathbf{y} - \mathbf{x}}{|\mathbf{x} - \mathbf{y}|^2}.$$

Let $\nu_y$ denote the unit outer normal at point $\mathbf{y} \in \Gamma$. The double layer potential with continuous density $\sigma$ is given by

$$K\sigma(\mathbf{x}) := \int_\Gamma \frac{\partial \Phi(\mathbf{x}, \mathbf{y})}{\partial \nu_y} \sigma(\mathbf{y}) ds_y = \int_0^1 \frac{\partial \Phi(\mathbf{x}, \mathbf{r}(t))}{\partial \nu_y} |\mathbf{r}'(t)| \sigma(\mathbf{r}(t)) dt, \quad \mathbf{x} \in \Omega, \tag{6.7}$$

where we assume $\Gamma$ is parametrized by $\mathbf{r}(t) : [0, 1] \to \mathbb{R}^2$.

Table 6.2: Numerical results for solving the Cauchy-like matrix when $\{x_k\}_{k=1}^n$ are distributed on three different curves.

| curve | $n$ | $\|u - \hat{u}\|/\|u\|$ | $\|Au - A\hat{u}\|/\|Au\|$ | $t_{\text{constr}}$ | $t_{\text{sol}}$ |
|---|---|---|---|---|---|
| [0, 1] | 1600 | $7.69 \times 10^{-12}$ | $5.56 \times 10^{-15}$ | 0.33 | 0.11 |
| | 3200 | $1.01 \times 10^{-09}$ | $3.87 \times 10^{-14}$ | 0.63 | 0.19 |
| | 6400 | $5.58 \times 10^{-11}$ | $5.58 \times 10^{-14}$ | 1.29 | 0.35 |
| | 12800 | $1.47 \times 10^{-08}$ | $5.87 \times 10^{-14}$ | 2.58 | 0.69 |
| honeybee | 1600 | $9.37 \times 10^{-11}$ | $8.09 \times 10^{-14}$ | 1.14 | 0.28 |
| | 3200 | $9.78 \times 10^{-10}$ | $5.60 \times 10^{-13}$ | 2.22 | 0.51 |
| | 6400 | $1.55 \times 10^{-09}$ | $9.16 \times 10^{-13}$ | 4.42 | 0.96 |
| | 12800 | $2.76 \times 10^{-09}$ | $1.54 \times 10^{-12}$ | 8.49 | 1.87 |
| snail | 1600 | $1.34 \times 10^{-11}$ | $1.02 \times 10^{-15}$ | 1.75 | 0.39 |
| | 3200 | $2.98 \times 10^{-11}$ | $6.61 \times 10^{-16}$ | 2.93 | 0.69 |
| | 6400 | $3.65 \times 10^{-10}$ | $5.88 \times 10^{-15}$ | 5.81 | 1.30 |
| | 12800 | $2.78 \times 10^{-10}$ | $6.27 \times 10^{-15}$ | 10.57 | 2.39 |

Given Dirichlet data $u_D \in C(\Gamma)$ in (6.6), we solve the following integral equation for $\sigma \in C(\Gamma)$:

$$(6.8) \qquad (K - \frac{1}{2}I)\sigma = u_D, \quad \text{on } \Gamma.$$

It is well-known ([39]) that the problem above for $\sigma \in C(\Gamma)$ is well-posed, and the corresponding double layer potential $u := K\sigma$ solves the interior Dirichlet problem (6.6).

Denote the kernel in the second integral in (6.7) by

$$(6.9) \qquad \kappa(s, t) := \frac{\partial \Phi(\mathbf{r}(s), \mathbf{r}(t))}{\partial \nu_y} |\mathbf{r}'(t)|.$$

It is easily seen that the smoothness of the above kernel function depends on the underlying boundary curve and it is known that ([2]) if $\mathbf{r} \in C^\infty([0, 1])$, then $\kappa \in C^\infty([0, 1] \times [0, 1])$.

Several Laplace problems (6.6) with the same exact solution but different domains are considered here. The first domain $\Omega$ is a ram head whose boundary curve $\Gamma$ is parametrized by $\mathbf{r}(t) = (r_1(t), r_2(t))$ for $t \in [0, 1]$:

$$(6.10) \qquad \begin{aligned} r_1(t) &= 2\cos(2\pi t), \\ r_2(t) &= 1 + \sin(2\pi t) - 1.4\cos^4(4\pi t). \end{aligned}$$

The second domain is a sunflower whose boundary curve $\Gamma$ is parametrized by:

$$(6.11) \qquad \begin{aligned} r_1(t) &= (1.3 + 1.25\cos(40\pi t))\cos(2\pi t), \\ r_2(t) &= (1.3 + 1.25\cos(40\pi t))\sin(2\pi t). \end{aligned}$$

We chose the Dirichlet data $u_D$ such that the exact solution of (6.6) is

$$u(\mathbf{x}) = \log|\mathbf{x} - \mathbf{x}_0|,$$

where the source point $\mathbf{x}_0 = (2, 1.5)$ is in the exterior of $\Omega$. Illustrations for the curves parametrized in (6.10) and (6.11) are shown in Fig. 6.3 and Fig. 6.4, respectively. We used Nyström method with trapezoidal rule to discretize (6.8). Since the curve $\Gamma$ and the kernel are both smooth, Nyström discretization converges with a convergence rate proportional to that of the quadrature rule.

As in Section 6.3, we applied the HSS matrix techniques to approximate and solve the resulting matrix from the Nyström discretization of the integral equation in (6.8). The adaptive partitioning based on bisection (cf. Remark A.1) was applied to a box covering the domain $\Omega$ in $\mathbb{R}^2$ and each box in leaf level contained no more than 50 quadrature points on $\Gamma$. Empty boxes were discarded during the partitioning. Illustrations of adaptive partitioning are shown in the right subfigures of Fig. 6.3 and Fig. 6.4, when 10240 quadrature points are in use. A binary tree $\mathcal{T}$ was then generated corresponding to each adaptive partitioning. In the construction of the HSS matrix, the bases for the farfield blocks were approximated by polynomial interpolation (with 25 interpolating points) with respective to the separation ratio 0.6 and the bases for the nearfield blocks are computed by the truncated SVD with the tolerance $\epsilon_{\text{svd}} = 10^{-11}$. In order to test the convergence of the discretization, for curves in Fig. 6.3, we compared the numerical solution $\hat{u}$ with the exact solution $u$ by evaluating them at point $\mathbf{x}^* = (0.1, 0.1)$ inside $\Omega$. For the curve in Fig. 6.4, the evaluation point is chosen as $\mathbf{x}^* = (1.5, 0)$ inside $\Omega$. The numerical results for the ram head and sunflower problems are shown in Table 6.3 and Table 6.4, respectively.

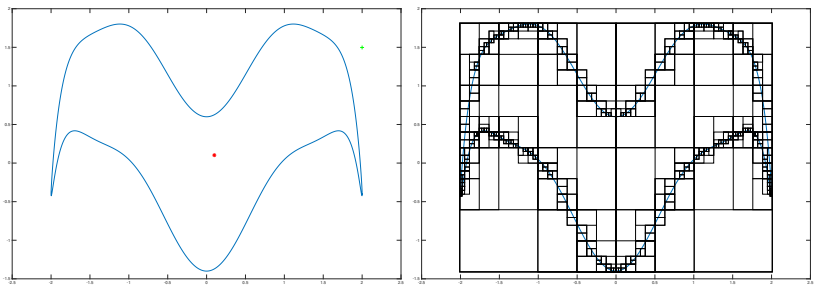

Fig. 6.3: Ram head domain for the Dirichlet problem (6.6) with source point and evaluation point marked as green '+' and red '*', respectively. Left: Original curve; Right: Adaptive partitioning of the ram head boundary curve for the case when 10240 quadrature points are used in Table 6.3.

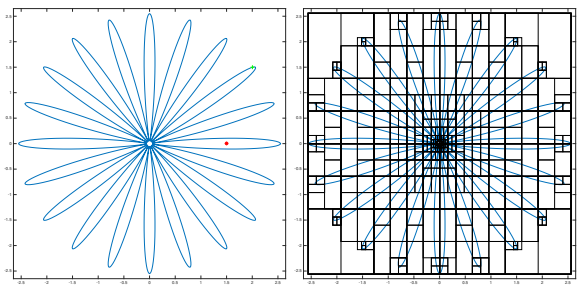

Fig. 6.4: Sunflower domain for the Dirichlet problem (6.6) with source point and evaluation point marked as green '+' and red '*', respectively. Left: Original curve; Right: Adaptive partitioning of the sunflower boundary curve for the case when 10240 quadrature points are used in Table 6.4.

From Table 6.3 and Table 6.4, it can be seen that the HSS matrix methods achieve linear complexity at both the construction and solution stages. In addition, it is worth noting that the numerical solutions for these Laplace Dirichlet problems converge exponentially fast regardless of the complicated geometries and the solver is quite robust. For example, in view of Table 6.4 which corresponds to the seemingly complicated geometry in Fig. 6.4, 10 digits of accuracy can be achieved using only 5120 quadrature points.

Table 6.3: Numerical results for solving a 2D Laplace Dirichlet problem in a ram head domain as shown in Fig. 6.3.

| $n$ | $|u(\mathbf{x}^*) - \hat{u}(\mathbf{x}^*)|$ | $\|A - \hat{A}\|_{\max}$ | cond(A) | $t_{\text{constr}}$ | $t_{\text{sol}}$ |
|---|---|---|---|---|---|
| 160 | $5.03 \times 10^{-08}$ | $1.06 \times 10^{-10}$ | $6.15 \times 10^{01}$ | 0.42 | 0.107 |
| 320 | $9.54 \times 10^{-11}$ | $6.81 \times 10^{-10}$ | $6.92 \times 10^{1}$ | 1.54 | 0.042 |
| 640 | $1.91 \times 10^{-12}$ | $7.98 \times 10^{-10}$ | $6.00 \times 10^{1}$ | 2.93 | 0.127 |
| 1280 | $8.22 \times 10^{-13}$ | $2.26 \times 10^{-09}$ | $6.02 \times 10^{1}$ | 5.41 | 0.103 |
| 2560 | $7.78 \times 10^{-13}$ | $3.90 \times 10^{-09}$ | $6.02 \times 10^{1}$ | 9.76 | 0.177 |
| 5120 | $1.50 \times 10^{-13}$ | $9.63 \times 10^{-09}$ | $6.02 \times 10^{1}$ | 18.74 | 0.282 |
| 10240 | $1.96 \times 10^{-12}$ | $1.01 \times 10^{-08}$ | $6.02 \times 10^{01}$ | 34.78 | 0.591 |

Table 6.4: Numerical results for solving a 2D Laplace Dirichlet problem in a sunflower domain as shown in Fig. 6.4.

| $n$ | $|u(\mathbf{x}^*) - \hat{u}(\mathbf{x}^*)|$ | $\|A - \hat{A}\|_{\max}$ | cond($A$) | $t_{\text{constr}}$ | $t_{\text{sol}}$ |
|---|---|---|---|---|---|
| 640 | $2.96 \times 10^{-02}$ | $1.20 \times 10^{-08}$ | $5.65 \times 10^{3}$ | 5.42 | 0.105 |
| 1280 | $1.25 \times 10^{-03}$ | $2.65 \times 10^{-08}$ | $2.23 \times 10^{3}$ | 16.80 | 0.256 |
| 2560 | $1.88 \times 10^{-06}$ | $3.56 \times 10^{-08}$ | $1.85 \times 10^{3}$ | 41.76 | 0.624 |
| 5120 | $1.02 \times 10^{-10}$ | $4.27 \times 10^{-07}$ | $1.72 \times 10^{4}$ | 102.65 | 1.429 |
| 10240 | $1.66 \times 10^{-11}$ | $3.80 \times 10^{-07}$ | $1.12 \times 10^{4}$ | 192.78 | 2.165 |
| 20480 | $8.03 \times 10^{-10}$ | $9.55 \times 10^{-07}$ | $6.86 \times 10^{3}$ | 316.68 | 3.219 |

**6.5. Nearly optimal compression.** In this section, we compare the exact numerical rank of the largest off-diagonal block of $A$ with the approximation rank obtained from SMASH. The numerical results show that the approximation rank is nearly optimal in the sense that it differs from the exact numerical rank by a small constant that is roughly independent of the kernel and the matrix size.

DEFINITION 6.1 ($\epsilon$-rank). *Let $\sigma_1 \geq \sigma_2 \geq \cdots \geq \sigma_r$ be singular values of a nonzero matrix $A$. Given a tolerance $\epsilon \in (0,1)$, the relative $\epsilon$-rank of a matrix $A$, denoted by $r_\epsilon(A)$, is the largest number $i$ such that $\sigma_i \geq \epsilon \sigma_1$.*

Consider the numerical examples in Section 6.4, where different curves give rise to different kernels according to (6.9). Let $i, j$ be two children of the root node. We focus on the (largest) off-diagonal block $A_{\mathbf{i}^{row} \times \mathbf{j}^{col}}$.

We consider three tolerances: $\epsilon = 10^{-3}, 10^{-6}, 10^{-10}$. We list the size of $A_{\mathbf{i}^{row} \times \mathbf{j}^{col}}$, the exact $\epsilon$-rank, and the approximation rank characterized by the size of $B_i$ generator with size($B_i$) := maximum between row size and column size. The results are reported in Table 6.5 .

Note that no a priori information is needed to determine the approximation rank as it is solely derived from the prescribed tolerance and the construction algorithm. Thus the numerical results also imply that the proposed method in Section 6.1 for choosing parameters is satisfactory.

**6.6. Storage cost.** In this section, we demonstrate the benefit of the special structure in the generators produced by SMASH. We store the HSS generators using the strategy mentioned at the end of Section 2.3. For comparison, we also compute the cost by storing the generators as dense matrices, denoted by HSS$_0$, as well as the storage cost for the original dense matrix. The test matrix

Table 6.5: Comparison of exact $\epsilon$-rank and approximation rank of $M = A_{\mathbf{i}^{row} \times \mathbf{j}^{col}}$ for the ram head geometry in Fig.6.3, the sunflower geometry in Fig.6.4 with $\epsilon = 10^{-3}, 10^{-6}, 10^{-10}$

| geometry | $n$ | size($M$) | $\epsilon = 10^{-3}$ | | $\epsilon = 10^{-6}$ | | $\epsilon = 10^{-10}$ | |
|---|---|---|---|---|---|---|---|---|
| | | | $r_\epsilon(M)$ | size($B_i$) | $r_\epsilon(M)$ | size($B_i$) | $r_\epsilon(M)$ | size($B_i$) |
| ram head | 1280 | 640 | 13 | 19 | 25 | 45 | 43 | 70 |
| | 2560 | 1280 | 13 | 18 | 25 | 45 | 43 | 71 |
| | 5120 | 2560 | 13 | 18 | 25 | 45 | 43 | 72 |
| | 10240 | 5120 | 13 | 19 | 25 | 45 | 43 | 72 |
| sun-flower | 1280 | 640 | 67 | 84 | 111 | 141 | 151 | 185 |
| | 2560 | 1280 | 83 | 117 | 159 | 187 | 213 | 251 |
| | 5120 | 2560 | 83 | 125 | 182 | 226 | 298 | 347 |
| | 10240 | 5120 | 83 | 123 | 187 | 237 | 328 | 380 |
| | 20480 | 10240 | 83 | 120 | 187 | 238 | 327 | 382 |

Table 6.6: A comparison of storage costs (MB) of HSS generators the storage strategy in Section 2.3 (SMASH) and the standard approach by storing dense generators (HSS$_0$) for a square matrix $A$ of order 10240.

| $\epsilon_{\text{far}}$ | $\epsilon_{\text{SVD}}$ | geometry | storage($A$) | HSS$_0$ | SMASH |
|---|---|---|---|---|---|
| $10^{-4}$ | $10^{-5}$ | ram head | 800 | 5.9 | 4.4 |
| $10^{-4}$ | $10^{-5}$ | sunflower | 800 | 35.2 | 11.0 |
| $10^{-10}$ | $10^{-11}$ | ram head | 800 | 23.2 | 8.8 |
| $10^{-10}$ | $10^{-11}$ | sunflower | 800 | 145.7 | 31.3 |

$A$ of order 10240 is derived from kernels in Section 6.4 and all entries are stored in double precision. The results are collected in Table 6.6 for different geometries and different approximation accuracy. The reduction in storage justifies the use of strong rank-revealing QR algorithm in the construction and we see that SMASH is quite cheap even when a high approximation accuracy is in use.

**7. Conclusion.** We presented a unified framework, called SMASH, to construct either an $n \times n$ HSS or $\mathcal{H}^2$ matrix with an $O(n)$ cost. One appealing feature of this scheme is its simple implementation which only requires a routine to compress far field blocks. In addition, SMASH can greatly reduce the memory cost relative to existing analytic construction schemes. The numerical experiments illustrated the efficiency and robustness of SMASH through a few examples with various point distributions and kernel matrices.

We plan to extend this scheme to highly oscillatory kernels and to develop approximate inverse-type preconditioners for solving the resulting linear systems with $\mathcal{H}^2$ matrix representations.

**Acknowledgments.** YX would like to thank Ming Gu for fruitful discussions about the strong rank revealing QR algorithm.

Fig. A.1: Illustration of the sets $\mathcal{N}_i$ used in HSS constructions.

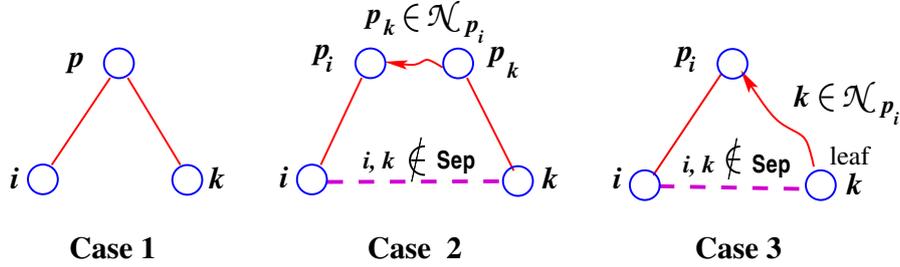

**Appendix A.** In building HSS structures, we need to collect nodes corresponding to blocks that are nearfield with respect to node $i$. The set of such nodes is denoted by $\mathcal{N}_i$ for node $i$. We set $\mathcal{N}_i = \emptyset$ when $i = root$. For the other cases, $\mathcal{N}_i$ is defined below where $p_i$ denotes the parent of $i$:

$$\text{(A.1)} \quad \mathcal{N}_i = \{k \in \mathcal{T} \text{ such that :} \quad \begin{array}{ll} \text{either} & k \text{ is a sibling of } i \\ \text{or} & p_k \in \mathcal{N}_{p_i} \text{ and } (i,k) \notin \text{Sep} \\ \text{or} & k \text{ is a leaf such that } k \in \mathcal{N}_{p_i} \text{ and } (i,k) \notin \text{Sep}\} \end{array}$$

Note that when $i$ is a child of root then $\mathcal{N}_{p_i}$ is empty and so only the first case can take place ($k$ is a sibling of $i$). We also remark that the third case ($k$ is a leaf such that $k \in \mathcal{N}_p$ and $(i,k) \notin \text{Sep}$) is required for non-uniform distributions and that it is empty if the distribution is uniform. It is easy to see that if $\Omega = [0,1]$ and $X = Y$ is uniformly distributed, $\mathcal{T}$ is a perfect binary tree. In addition, if the separation ratio is set to $\tau = 0.5$, then for any nonroot node $i$, $\mathcal{N}_i$ contains at most two nodes.

For each node $i$, let $\mathbf{i}$ denote the index set of the points in $X \cap \Omega_i$. Similarly, $\mathbf{j}$ represents the index set of the points in $Y \cap \Omega_j$. Namely, $X_{\mathbf{i}} = X \cap \Omega_i$ and $Y_{\mathbf{j}} = Y \cap \Omega_j$.

REMARK A.1. *Since the HSS structure [21, 22] is associated with a binary tree regardless of the dimension of the problem (see Section 3.2), to construct HSS matrices, bisection is used throughout the adaptive partitioning procedure. For example, given a domain or a curve enclosed in a square in* $\mathbb{R}^2$, *we use bisection in the horizontal direction and the vertical direction alternatively in consecutive stages of the adaptive partitioning, i.e., if horizontal bisection is used at partitioning level $l$, then vertical bisection will be employed at level $l+1$. The numerical experiments in Section 6.3 and Section 6.4 provide illustrations. This partitioning strategy corresponds to the* geometrically regular clustering (cf.[12]), *and can be generalized into the* geometrically balanced clustering (cf.[12]).